\documentclass[leqno,12pt]{amsart}
\usepackage{amssymb,amsmath}

\calclayout

\newtheorem{theorem}{Theorem}
\newtheorem{corollary}[theorem]{Corollary}

\newtheorem{lemma}[theorem]{Lemma}

\newcommand{\bn}{{\mathbf n}}
\newcommand{\bX}{{\mathbf X}}

\newcommand{\cI}{{\mathcal I}}
\newcommand{\cL}{{\mathcal L}}
\newcommand{\cO}{{\mathcal O}}
\newcommand{\cV}{{\mathcal V}}
\newcommand{\cX}{{\mathcal X}}

\newcommand{\mN}{{\mathbb N}}
\newcommand{\mP}{{\mathbb P}}
\newcommand{\mZ}{{\mathbb Z}}

\newcommand{\oT}{\overline{T}}
\newcommand{\oB}{\overline{B}}

\newcommand{\tB}{{\tilde B}}
\newcommand{\tG}{{\tilde G}}
\newcommand{\tP}{{\tilde P}} 
\newcommand{\tT}{{\tilde T}}
\newcommand{\tcX}{\tilde{\mathcal X}}

\def\deuxind#1#2{{\buildrel{\scriptstyle #1}\over{#2}}}

\title{Large Schubert varieties}
\author{Michel~Brion}
\address{Universit\'e de Grenoble I\\
D\'epartement de Math\'ematiques\\
Institut Fourier, UMR 5582 du CNRS\\
38402 Saint-Martin d'H\`eres Cedex, France}
\email{Michel.Brion@ujf-grenoble.fr}
\author{Patrick~Polo}
\address{Universit\'e Paris Nord\\
D\'epartement de Math\'ematiques\\
L.A.G.A., UMR 7539 du CNRS\\
93430 Villetaneuse, France}
\email{polo@math.univ-paris13.fr}
\date{}
\begin{document}

\begin{abstract}
For a semisimple adjoint algebraic group $G$ and a Borel subgroup $B$,
consider the double classes $BwB$ in $G$ and their closures in the
canonical compactification of $G$: we call these closures large
Schubert varieties. We  show that these varieties are normal and
Cohen-Macaulay; we describe their Picard group and the spaces of
sections of their line bundles. As an application, we construct
geometrically van der Kallen's filtration of the algebra of regular
functions on $B$. We also construct a degeneration of the flag variety
$G/B$ embedded diagonally in $G/B\times G/B$, into a union of Schubert
varieties. This leads to formulae for the class of the diagonal in 
$T$-equivariant $K$-theory of $G/B\times G/B$, where $T$ is a maximal
torus of $B$.
\end{abstract}

\maketitle

\section*{Introduction}

Consider an adjoint semisimple algebraic group $G$ and a Borel
subgroup $B$. The Schubert varieties are the images in $G/B$ of the
closures in $G$ of double classes $BwB$. These varieties are generally
singular, but all of them are normal and Cohen-Macaulay \cite{R1}. The 
spaces of sections of line bundles over Schubert varieties play an
important r\^ole in representation theory, see for example \cite{J},
\cite{vdK2}.
\smallskip

The group $G$ has a canonical smooth $G\times G$-equivariant
completion $\bX$, constructed by DeConcini and Procesi \cite{DP} in
characteristic zero, and by Strickland \cite{S} in arbitrary
characteristics. In this paper, we study the closures of double
classes $BwB$ in $\bX$; we call them {\sl large Schubert varieties}.
\smallskip

These varieties are highly singular: by \cite[2.2]{B}, their singular
locus has codimension two, apart from trivial exceptions.
However, we show that large Schubert varieties are normal and
Cohen-Macaulay (Corollary \ref{sur} and Theorem \ref{cm1}). Further,
their Picard group is isomorphic to the weight lattice (Theorem
\ref{pic}).
\smallskip

Large Schubert varieties have an obvious relation to usual Schubert
varieties: the latter are quotients by $B$ of an open subset of the
former. A more hidden connexion arises by intersecting a large
Schubert variety $X$ with the unique closed $G\times G$-orbit $Y$ in
$\bX$. As $Y$ is isomorphic to $G/B\times G/B$ by \cite{S}, $X\cap Y$
is a union of Schubert varieties in $G/B\times G/B$.
\smallskip

The space $X\cap Y$ is generally reducible; its irreducible components
were described in \cite{B}, e.g. those of $\overline{B}\cap Y$ are 
parametrized by the Weyl group. We show that the 
scheme-theoretic intersection $X\cap Y$ is reduced and Cohen-Macaulay
(Corollaries \ref{reg}, \ref{cm2}). Together with a construction of
\cite{B}, this leads to a degeneration of the diagonal in 
$G/B\times G/B$ into a union of Schubert varieties, and then to
formulae for the class of the diagonal in equivariant $K$-theory
(Theorem \ref{flat} and Corollary \ref{diag}).
\smallskip

Let $\tB$ be the preimage of $B$ in the simply-connected cover $\tG$ 
of $G$. Then the space of sections of each line bundle over a large
Schubert variety $X$ is a $\tB\times\tB$-module, endowed with a
natural filtration by order of vanishing of sections along $Y$. We
decompose the associated graded into a direct sum of spaces of
sections of line bundles over $X\cap Y$; the latter
$\tB\times\tB$-modules are indecomposable (Theorem \ref{fil} 
and Corollary \ref{ind}). In the case where $X=\oB$, these
modules can be seen as degenerations of induced $\tG$-modules, see
Corollary \ref{sep}.
\smallskip

As a consequence, we recover van der Kallen's filtration of the
affine algebra of $\tB$ \cite{vdK1} in a geometric way
(Theorem \ref{vdk}). 
For this, consider regular functions on $B$ as rational functions on
its closure $\oB$ with poles along the boundary. The factors of the
filtration by order of poles are spaces of sections of line bundles on
$\oB\cap Y$. In particular, as a $B\times B$-module, the affine
algebra of $B$ admits a Schubert filtration in the sense of
\cite{P}. Filtering further by ordering the irreducible components of
$\oB\cap Y$ gives back the filtration of \cite[1.13]{vdK1} for
$B$. This generalizes to regular functions over $\tB$, by decomposing
them into sums of sections of line bundles over $B$.
\smallskip

Our proofs rely on the method of Frobenius splitting:
following the approach of \cite{LT}, we show that $\bX$ is
Frobenius split compatibly with all large Schubert varieties
(Theorem \ref{split}). The normality of large Schubert varieties
is a direct consequence of this fact: it is easy to see that they are
smooth in codimension one, and that their depth is at least two 
(a regular sequence being provided by the ``boundary divisors'' of
$\bX$.) 
\smallskip

The proof that large Schubert varieties are Cohen-Macaulay is
much more involved. As for usual Schubert varieties \cite{R1}, we
proceed by ascending induction on the dimension; but here the argument 
begins with $\oB$ (instead of the point) which is handled through its
intersection with $Y$. It would be interesting to obtain an
equivariant desingularization of $\oB$; then the classical
construction of Bott-Samelson-Demazure would give equivariant
resolutions of all large Schubert varieties. The present work raises
many other questions, e.g. is there a standard monomial theory
for large Schubert varieties ?
\smallskip

One may also ask for extensions of our results to orbit closures of
Borel subgroups in complete symmetric varieties \cite{DP}, or, more
generally, in regular embeddings of spherical homogeneous spaces
\cite{B}. It turns out that closures of orbits of maximal rank are
normal, and that their intersection with any irreducible component of
the boundary is reduced; further, these intersections can be described
in terms of the Weyl group. But many orbit closures of smaller rank
are neither normal, nor Cohen-Macaulay (see \cite{K2} for the notion
of rank; all large Schubert varieties have maximal rank.) This will be
developed elsewhere.

\section{The canonical completion of a semisimple adjoint group}

We begin by introducing notation and recalling some properties
of group completions.

Let $G$ be a connected adjoint semisimple algebraic group over an
algebraically closed field $k$ of arbitrary characteristic. Let $B$
and $B^-$ be opposite Borel subgroups of $G$, with common torus
$T$. Let $\cX$ be the character group of $T$; we identify 
$\cX$ with the character groups of $B$ and $B^-$. Let $W$ be the
Weyl group of $T$, and let $\Phi$ be the root system of $(G,T)$ with
subsets of positive (resp. negative) roots $\Phi^+$, $\Phi^-$ defined
by $B$, $B^-$. Let $\Delta=\{\alpha_1,\ldots,\alpha_r\}$ be the set of
simple roots, where $r$ is the rank of $G$; let $s_1,\ldots, s_r\in W$
be the simple reflections. The corresponding length function on $W$ is
denoted by $\ell$. Let $w_0$ be the element of maximal length in
$W$. Set $N:=\ell(w_0)$, the number of positive roots.

We denote by $\tG$ the simply connected covering of $G$, and by
$\tB$, $\tT$, $\ldots$ the preimages of $B$, $T$, $\ldots$
in $\tG$. The character group of $\tT$ is denoted by
$\tcX$; it is the weight lattice of $\Phi$ with basis the
set of fundamental weights $\omega_1,\ldots,\omega_r$. 
The monoid generated by these weights is the set $\tcX^+$ of dominant
weights. Let $\leq$ denote the usual partial order on $\tcX$ defined by:
$\lambda\leq\mu$ if there exist non-negative integers $n_1,\ldots,n_r$
such that $\mu-\lambda= n_1\alpha_1+\cdots+n_r\alpha_r$.

\smallskip 
By \cite{DP} and \cite{S}, $G$ admits a completion $\bX$ satisfying
the following properties:

\noindent
(i) $\bX$ is a smooth projective variety, and the action of
$G\times G$ on $G$ by left and right multiplication extends to
$\bX$.

\noindent
(ii) The boundary $\bX - G$ is a union of $r$ smooth
irreducible divisors $D_1$, $\ldots$, $D_r$ with normal crossings.

\noindent
(iii) Each $G\times G$-orbit closure in $\bX$ is the
transversal intersection of the boundary divisors which contain it.

\noindent
(iv) The intersection $D_1\cap\cdots\cap D_r := Y$ is the unique
closed $G\times G$-orbit in $\bX$; it is isomorphic to 
$G/B\times G/B$.

Further, any completion of $G$ satisfying (i), (ii) and (iii)
dominates $\bX$, and any normal completion of $G$ with a unique closed 
orbit is dominated by $\bX$ (this follows from embedding theory of
homogeneous spaces, see \cite{K1}). We will call $\bX$ the 
{\it canonical completion} of $G$.

For $w\in W$, consider the double class $BwB$ in $G$, and its closure
in $\bX$. We denote this closure by $\bX(w)$, and we call it a 
{\it large Schubert variety}. On the other hand, we denote by $S(w)$
the usual Schubert variety, that is, the closure in $G/B$ of $BwB/B$. 
In other words, $S(w)$ is the image in $G/B$ of the intersection
$\bX(w)\cap G$.

The intersections of large Schubert varieties with 
$G\times G$-orbits were studied in \cite[\S 2]{B}. In particular, we
have the following decomposition of 
$$
Z(w):=\bX(w)\cap Y
$$
into irreducible components (which must be Schubert varieties in
$G/B\times G/B$): 
$$
Z(w) = \bigcup_{\deuxind{x\in W}{\ell(wx)=\ell(w)+\ell(x)}}
S(wx)\times S(xw_0).
$$
For $w=1$ (the identity of $W$), we denote $\bX(w)$ by $\oB$, and
$Z(w)$ by $Z$. Then
$$
Z=\bigcup_{x\in W} S(x)\times S(xw_0).
$$

The large Schubert varieties of codimension one in $\bX$ are
$\bX(w_0s_1),\ldots,\bX(w_0s_r)$. They are the irreducible 
$B\times B$-stable divisors in $\bX$ which are not $G\times G$-stable,
or, equivalently, which do not contain $Y$. By
\cite[Proposition 4.4]{DS}, the divisor class group of $\bX$ is freely
generated by the classes of $\bX(w_0s_1),\ldots,\bX(w_0s_r)$. 
On the other hand, the line bundles on $\bX$ are described in 
\cite[\S 2]{S} (see also \cite[\S 4]{DS}). We now recall this
description; a generalization to large Schubert varieties will be
obtained in Section 3. 

For $\lambda$ and $\mu$ in $\tcX$, let $\cL_Y(\lambda,\mu)$ be the
corresponding homogeneous line bundle on $Y=G/B\times G/B$. The map 
$(\lambda,\mu)\mapsto\cL_Y(\lambda,\mu)$ identifies the Picard group
${\rm Pic}\;Y$ with $\tcX\times\tcX$. Now the
restriction $res_Y:{\rm Pic}\;\bX\to{\rm Pic}\;Y$ is
injective, and its image consists of the $\cL_Y(\lambda,-w_0\lambda)$, 
for $\lambda\in\tcX$.

We denote by $\cL_{\bX}(\lambda)$ the line bundle on $\bX$ such that 
$res_Y\cL_{\bX}(\lambda)=\cL_Y(\lambda,-w_0\lambda)$. This identifies
${\rm Pic}\;\bX$ with $\tcX$; we can index the boundary divisors
$D_1,\ldots, D_r$ so that the classes of the corresponding line
bundles are $\cL_{\bX}(\alpha_1),\ldots,\cL_{\bX}(\alpha_r)$. Then
each $\cL_{\bX}(\alpha_i)=\cO_{\bX}(D_i)$ has a section $\sigma_i$ with
divisor $D_i$; this section is unique up to scalar multiplication.

Because $\tG$ is semisimple and simply connected, each line bundle
$\cL_{\bX}(\lambda)$ has a unique $\tG\times\tG$-linearization. Thus,
each space $H^0(\bX,\cL_{\bX}(\lambda))$ is a $\tG\times\tG$-module,
which we denote by $H^0(\bX,\lambda)$. Similarly, we denote
$H^0(Y,\cL_Y(\lambda,\mu))$ by $H^0(Y,\lambda,\mu)$. This
$\tG\times\tG$-module 
is isomorphic to $H^0(G/B,\lambda)\boxtimes H^0(G/B,\mu)$
where $\boxtimes$ denotes the external tensor product.  

Observe that the section $\sigma_i$ of $\cL_{\bX}(\alpha_i)$ is
$\tG\times\tG$-invariant. This is the starting point for an analysis
of the $\tG\times\tG$-module $H^0(\bX,\lambda)$ for arbitrary
$\lambda$, see \cite[\S 2]{S}; the results will be generalized to
large Schubert varieties in Section 3. Here we will need the following 

\begin{lemma}\label{div}
For any dominant weight $\lambda$, the line bundle
$\cL_{\bX}(\lambda)$ has a global section $\tau_{\lambda}$,
eigenvector of $\tB\times\tB$ of weight $(-w_0\lambda,\lambda)$. This
section is unique up to scalar, and its divisor is
$$
\sum_{i=1}^r \langle\lambda,\alpha_i^{\vee}\rangle
\bX(w_0s_i).
$$
\end{lemma}

\begin{proof}
By \cite[\S 4]{DS}, the line bundle on $\bX$ associated with the
divisor $\bX(w_0s_i)$ is $\cL_{\bX}(\omega_i)$. Let $\tau_i$ be the
canonical section of this line bundle; then $\tau_i$ is an
eigenvector of $\tB\times\tB$, because its divisor $\bX(w_0s_i)$ is 
$B\times B$-stable. The closure in $\tG$ of $\tB w_0s_i\tB$ is the 
divisor of a regular function on $\tG$, eigenvector of $\tB\times\tB$ of
weight $(-w_0\omega_i,\omega_i)$, and unique up to scalar
multiplication. Thus, the weight of $\tau_i$ is
$(-w_0\omega_i,\omega_i)$. As the $\cL_{\bX}(\omega_i)$ generate the
Picard group of $\bX$, the existence of $\tau_{\lambda}$ and the
formula for its divisor follow immediately. Finally, uniqueness of
$\tau_{\lambda}$ up to scalar is a consequence of the fact that
$B\times B$ has a dense orbit in $\bX$.
\end{proof}

\section{Compatible Frobenius splitting and applications}

In the beginning of this section, we assume that $k$ has
characteristic $p>0$. 
For a scheme $X$ over $k$, we denote by $F:X\to X$ the absolute
Frobenius morphism. Recall that $X$ is {\it Frobenius split} if the
map $F:\cO_X\to F_*\cO_X$ is split, that is, if there exists
$\sigma\in {\rm Hom}_{\cO_X}(F_*\cO_X,\cO_X)$ such that 
$\sigma\circ F$ is the identity. Let $Y\subseteq X$ be a closed
subscheme with ideal sheaf $\cI_Y$; then a splitting $\sigma$ is
{\it compatible with} $Y$ if $\sigma(F_*\cI_Y)$ is contained in $\cI_Y$.

By \cite[\S 3]{S}, the canonical completion $\bX$ is Frobenius split
compatibly with all $G\times G$-orbit closures. We will need the
following refinement of this result.

\begin{theorem}\label{split}
$\bX$ is Frobenius split compatibly with all $G\times G$-orbit
closures and all subvarieties $\bX(w)$ and $(w_0,w_0)\bX(w)$, for $w\in W$.
\end{theorem}

\begin{proof}
Let $St=H^0(G/B,(p-1)\rho)$ be the Steinberg module for $\tG$; it is a
simple, self-dual $\tG$-module \cite[II.2.5, II.3.18]{J}. On the other
hand, the line bundle $\cL_{\bX}((p-1)\rho)$ is 
$\tG\times\tG$-linearized by construction of $\bX$, and the
$\tG\times\tG$-module $H^0(\bX,(p-1)\rho)$ 
contains an eigenvector of $\tB\times\tB$ of weight $(p-1)(\rho,\rho)$,
unique up to scalar, by Lemma \ref{div}. Further, the image of this
eigenvector under restriction to $Y$ is non-zero, since no $\bX(w_0s_i)$
contains $Y$. Using Frobenius reciprocity \cite[I.3.4]{J} and
self-duality of $St$, we obtain a $\tG\times\tG$-homomorphism 
$$
f:St\boxtimes St \to H^0(\bX,(p-1)\rho)
$$
such that the composition
$$
res_Y\circ f:St\boxtimes St\to H^0(Y,(p-1)(\rho,\rho))
$$
is non-zero. Since the $\tG\times\tG$-module
$H^0(Y,(p-1)(\rho,\rho))$ is isomorphic to $St\boxtimes St$, hence
simple, it follows that $res_Y\circ f$ is an isomorphism.

We thus obtain a $\tG\times \tG$-homomorphism
$$
f^2:(St\boxtimes St)^{\otimes 2}\to
H^0(\bX,2(p-1)\rho),
$$
$$
x_1\boxtimes y_1\otimes x_2\boxtimes y_2\mapsto
f(x_1\boxtimes y_1)f(x_2\boxtimes y_2).
$$
Moreover, the composition
$$
res_Y\circ f^2:(St\boxtimes St)^{\otimes 2}\to
H^0(Y,2(p-1)(\rho,\rho))
$$
is surjective, because the product map
$$
H^0(Y,(p-1)(\rho,\rho))^{\otimes 2}\to H^0(Y,2(p-1)(\rho,\rho))
$$
is \cite[II.14.20]{J}. Now, by \cite[2.1,2.3]{LT}, there is a natural
$\tG\times \tG$-isomorphism 
$$
{\rm Hom}(F_*\cO_Y,\cO_Y) \buildrel\cong\over\longrightarrow
H^0(Y,2(p-1)(\rho,\rho))
$$
and there is a unique $\tG\times\tG$-homomorphism (up to a constant)
$$
\varphi:(St\boxtimes St)^{\otimes 2}\to
{\rm Hom}(F_*\cO_Y,\cO_Y).
$$
Further, for $a$ and $b$ in $St\boxtimes St$, the map 
$\varphi(a\otimes b)$ is a splitting of $Y$ (up to a constant) if and
only if  $\langle a,b\rangle\neq 0$ where
$\langle,\rangle$ is the $\tG\times\tG$-invariant bilinear form on
$St\boxtimes St$. Finally, if $a=s^{p-1}$ and $b=t^{p-1}$ for sections 
$s$, $t$ of $\cL_Y(\rho,\rho)$, then the zero subschemes $Z(s)$,
$Z(t)$ in $Y$ are compatibly $\varphi(a\otimes b)$-split.

Because $res_Y\circ f^2$ is a surjective
$\tG\times \tG$-homomorphism, we can identify it with $\varphi$. Let
$v_+$ (resp. $v_-$) be a highest (resp. lowest) weight vector in
$H^0(G/B,\rho)$. Set $s:=v_+\boxtimes v_+$, $t:=v_-\boxtimes v_-$,
$a:=s^{p-1}$ and $b:=t^{p-1}$. Then $a$, $b$ are in $St\boxtimes St$
and they satisfy $\langle a,b\rangle\neq 0$. Thus,  
$res_Y\circ f^2(a\otimes b)$ splits $Y$ compatibly with $Z(s)$
and $Z(t)$.

Set $\tau:=\varphi(a\otimes b)$ and consider
$$
\sigma:=\tau\prod_{i=1}^r \sigma_i^{p-1},
$$
a global section of
$\cL_{\bX}((p-1)(2\rho+\sum_{i=1}^r\alpha_i))$. Recall from 
\cite[\S3]{S} that the
canonical sheaf of $\bX$ is 
$$
\omega_\bX=\cL_{\bX}(-2\rho-\sum_{i=1}^r\alpha_i).
$$
Thus,
$\sigma\in\Gamma(\bX,\omega_\bX^{1-p})
\simeq{\rm Hom}_{\cO_{\bX}}(F_*\cO_{\bX},\cO_{\bX})$.
By \cite[Th. 3.1]{S}, $\sigma$ splits $\bX$
compatibly with $D_1,\ldots,D_r$.

Set $\tau_+=f(v_+^{p-1}\boxtimes v_+^{p-1})$ and
$\tau_-=f(v_-^{p-1}\boxtimes v_-^{p-1})$. Then $\tau_+$ and
$\tau_-$ are in $H^0(\bX,(p-1)\rho)$, and $\tau_+$ (resp. 
$\tau_-$) is an eigenvector of $\tB\times\tB$
(resp. $\tB^-\times\tB^-$) of weight $(p-1)(\rho,\rho)$
(resp. $-(p-1)(\rho,\rho)$). By Lemma \ref{div}, we have
$$
{\rm div}(\tau_{\pm})=(p-1)\bX_{\pm}
$$
where $\bX_+$ is the sum of the classes of the $\bX(w_0s_i)$
over all simple reflections $s_i$, and
$\bX_-=(w_0,w_0)\bX_+$. Thus, $\sigma$ splits $\bX$ compatibly with
$\bX_+$ and $\bX_-$.
This implies the theorem, as in the proof of
\cite[Th. 3.5(i)]{R1}. Namely, one uses \cite[Lemma 1.11]{R1} and the
fact that each $\bX(w)$ is an irreducible component of an iterated
intersection of irreducible components of $\bX_+$.
\end{proof}

\begin{corollary}\label{sur} Let $char(k)$ be arbitrary. 

\noindent (i) For any dominant weight $\lambda$ and for any
intersection $X$ of large Schubert varieties and of boundary divisors,
the restriction map
$$
res_X:H^0(\bX,\lambda)\to H^0(X,\lambda)
$$
is surjective. Further, $H^i(X,\lambda)=0$ for $i\geq 1$.

\smallskip\noindent (ii) Any intersection of large
Schubert varieties and of boundary divisors is reduced.
\end{corollary}

\begin{proof} 
Let us prove (i). Since the divisor $\bX_+$ is ample, 
this is a consequence of \cite[Proposition 1.13(ii)]{R2} when 
$char(k) = p > 0$. Moreover, since $G$, $B$ are defined over $\mZ$, it 
follows from the construction of $\bX$ (\cite{S}) that $\bX$, the
boundary divisors $D_i$ and the large Schubert varieties $\bX(w)$ are
all defined and flat over some open subset of 
${\rm Spec}\;\mZ$ (in fact, they are defined over $\mZ[1/2]$ by
\cite{DS}.) Therefore, by the semicontinuity theorem, (i)  
holds in characteristic zero as well. 

Moreover, by the proof of \cite[Th. 3]{R1}, (ii) follows (in
arbitrary characteristic) from Theorem \ref{split}.   
\end{proof}  

For $1\leq i\leq r$, multiplication by $\sigma_i$ (a section of
$\cL_{\bX}(\alpha_i)$ with divisor $D_i$) defines an exact sequence
$$
0\to\cO_{\bX}(-D_i)\to\cO_{\bX}\to\cO_{D_i}\to 0.
$$
Because $Y$ is the transversal intersection of $D_1,\ldots, D_r$, the
image of the map
$$
(\sigma_1,\ldots,\sigma_r):\bigoplus_{i=1}^r\cO_{\bX}(-D_i)\to\cO_{\bX}
$$
is the ideal sheaf $\cI_Y$.

\begin{corollary}\label{reg} Again, let $char(k)$ be arbitrary. 

\noindent
(i) The ideal sheaf of the set-theoretic intersection
$Z(w)=\bX(w)\cap Y$ in $\bX(w)$ is generated by the image of
$(\sigma_1,\ldots,\sigma_r)$. 

\noindent
(ii) $\sigma_1,\ldots,\sigma_r$ form a regular sequence in $\bX(w)$.

\noindent
(iii) $\bX(w)$ is normal.
\end{corollary}

\begin{proof}  
By Corollary \ref{sur}, the scheme-theoretic intersection 
$\bX(w)\cap Y$ is reduced; this is equivalent to (i).

For (ii), we have to check that the image of each $\sigma_j$ in
$\cO_{\bX(w)}/(\sigma_1,\ldots,\sigma_{j-1})$ is not a
zero divisor. But the scheme-theoretic intersection
$$
\bX(w)\cap\bigcap_{i=1}^{j-1} D_i:=\bX(w)_{<j}
$$
is reduced. Further, by \cite[Th. 2.1]{B}, each irreducible
component of $\bX(w)_{<j}$ has codimension $j-1$ in $\bX(w)$, and is
not contained in $D_j$. Thus, the restriction of $\sigma_j$ to
$\bX(w)_{<j}$ does not vanish on any such component. It follows that
$\sigma_j$ is not a zero divisor in
$\cO_{\bX(w)_{<j}}=\cO_{\bX(w)}/(\sigma_1,\ldots,\sigma_{j-1})$.

For (iii), observe that $\bX(w)\cap G$ is smooth in codimension
one, as the preimage in $G$ of a Schubert variety in
$G/B$ (this goes back to Chevalley \cite[Cor., p. 10]{C}.)
Further, the intersection $\bX(w)\cap D_i$ is reduced for 
$1\leq i\leq r$, so that each irreducible component of this
intersection contains smooth points of $\bX(w)$. Thus, $\bX(w)$ is
smooth in codimension one (this also follows from 
\cite[Cor. 2.1]{B}.) By Serre's criterion, it is enough to prove that
$\bX(w)$ has depth at least two.

Because $B\times B$ acts on $\bX(w)$ with finitely many orbits and a
unique fixed point $y$ (the base point of $Y=G/B\times G/B$), it
suffices to prove that ${\bf X}(w)$ has depth at least two at
$y$. This is clear if $r\geq 2$, because the local equations of
$D_1,\ldots,D_r$ at $y$ form a regular sequence in the local ring
$\cO_{\bX(w),y}$. On the other hand, if $r=1$ then each $\bX(w)$ is
smooth. We have indeed $\tG={\rm SL}(2)$, $G={\rm PGL}(2)$, and $X$
is the projectivization of the space of $2\times 2$ matrices where $G$
acts by left and right multiplication. So $\bX(w)$ is either $\bX$ or
the projectivization of the subspace of upper triangular matrices.
\end{proof}

\section{Line bundles on large Schubert varieties}\label{section-pic} 

In this section, we describe the Picard group of large Schubert
varieties, and the spaces of global sections of line bundles on these
varieties.

\begin{theorem}\label{pic}
For any $w\in W$, the restriction map 
$$
res_{\bX(w)}:{\rm Pic}\;\bX\to{\rm Pic}\;\bX(w)
$$ 
is bijective. Further, the line bundle
$\cL_{\bX(w)}(\lambda)$ is generated by its global sections
$($resp. ample$)$ if and only if $\lambda$ is dominant
$($resp. dominant regular$)$.
\end{theorem}

\noindent
{\sc Remarks.} 1) We will see in Corollary \ref{effective} that
$\cL_{\bX(w)}(\lambda)$ admits nontrivial global sections if and only
if $\lambda$ is in the monoid generated by all simple roots and
fundamental weights.

\noindent
2) It is proved in \cite[\S 2]{S} that $\cL_{\bX}(2\lambda)$ is very
ample for any regular dominant $\lambda$. In fact, one can check that
$\cL_{\bX}(\lambda)$ is already very ample, using Corollary \ref{sur}.

\smallskip

\begin{proof}  
We will use the duality between line bundles and curves:
each closed curve $C$ in $\bX(w)$ defines an additive map 
${\rm Pic}\;\bX(w)\to\mZ, L\mapsto (L\cdot C)$ where $(L\cdot C)$ is
the degree of the restriction of $L$ to $C$. In fact, $(L\cdot C)$
only depends on the classes of $L$ and $C$ up to rational
equivalence. Further, $C$ is rationally equivalent to a positive
integral combination of closed irreducible $B\times B$-stable curves
\cite{FMSS}.

Examples of such curves are the ``Schubert curves'' 
$C(\alpha_i):=S(s_i)\times S(1)$ and 
$C'(\alpha_i):=S(1)\times S(s_i)$ in $G/B\times G/B$. 
Note that 
$$
(\cL_{\bX(w)}(\lambda)\cdot C(\alpha_i)) =
(\cL_{\bX(w)}(\lambda)\cdot C'(\alpha_i))=
\langle\lambda,\alpha_i^\vee\rangle
$$ 
for all $\lambda\in\tcX$. We first show the following

\begin{lemma}\label{curves}
The closed irreducible $B\times B$-stable curves in $\bX$ are the
$C(\alpha_i)$ and $C'(\alpha_i)$  for $1\leq i\leq r$. They are
contained in $\oB$. Further, each $C(\alpha_i)$ is rationally
equivalent in $\oB$ to $C'(-w_0\alpha_i)$.
\end{lemma}

\begin{proof}
The first assertion follows from the description of all 
$B\times B$-orbits in $\bX$ given in \cite[2.1]{B}. And as 
$\oB\cap Y$ contains both $S(w_0)\times S(1)$ and 
$S(1)\times S(w_0)$, it contains the $C(\alpha_i)$ and
$C'(\alpha_i)$. 

For the latter assertion, we begin by the case where 
$G={\rm PGL}(2)$. Then we saw that $\bX=\mP^3$ and
$\oB=\mP^2$. Further, $Y$ is a smooth quadric in $\mP^3$, and
both $C(\alpha)$, $C'(-w_0\alpha)$ are embedded lines. Thus, they
are rationally equivalent in $\mP^2$. 

The general case reduces to the previous one, as follows. Set
$X_i:=\cap_{j\neq i}D_j$, then $X_i$ is the closure of a unique
$G\times G$-orbit $X_i^0$ in $X$. Let $P_i$ be the parabolic subgroup
generated by $B$ and $s_i$; let $Q_i$ be the opposite parabolic
subgroup containing $B^-$, and let $L_i$ be their common Levi
subgroup. Then the $G\times G$-variety $X_i$ fibers equivariantly over
$G/P_i\times G/Q_i$, with fiber the canonical completion of the
adjoint group $L_i/Z(L_i)$ (this follows e.g. from 
\cite[Th. 3.16]{DS}). This group is isomorphic to ${\rm PGL}(2)$. 
Set $w_0Q_iw_0:=P_j$, the parabolic subgroup generated by $B$ and
$w_0s_iw_0$. Now $X_i$ fibers equivariantly over $G/P_i\times G/P_j$
and the fiber over the base point is a closed $B\times B$-stable
subvariety $F_i$ of $X_i$, isomorphic to $\mP^3$. Restricting this
fibration to $Y\subset X_i$, we obtain the canonical map
$G/B\times G/B\to G/P_i\times G/P_j$. Thus, $F_i$ contains both
$C(\alpha_i)=P_i/B \times B/B$ and $C(-w_0\alpha_i)=B/B\times P_j/B$.
Further, $B\times B$ has a unique closed orbit $\cO_i$ in
$X_i^0$; and $\cO_i$ is contained in $F_i\cap\oB$ (because $\oB$ meets 
all $G\times G$-orbits). Thus, the closure of $\cO_i$ in $X_i$ is
isomorphic to $\mP^2$, and contains both $C(\alpha_i)$ and
$C'(-w_0\alpha_i)$ as embedded lines. 
\end{proof}

We return to the proof of Theorem \ref{pic}.
For injectivity, let $\lambda$ be a weight such that the restriction of
$\cL_{\bX}(\lambda)$ to $\bX(w)$ is trivial. Then the restriction of
$\cL_{\bX}(\lambda)$ to each $C(\alpha_i)$ is trivial. It follows that
$\langle\lambda,\alpha_i^\vee\rangle = 0$ for $1\leq i\leq r$,
and that $\lambda=0$. 

For surjectivity, we first prove that the abelian group 
${\rm Pic}\;\bX(w)$ is free of finite rank. For this, we identify
${\rm Pic}\;\bX(w)$ to the group of all Cartier divisors on $\bX(w)$
up to rational equivalence (this holds because $\bX(w)$ is normal.)
Let $y$ be the $B\times B$-fixed point of $Y$. Let $\bX_y$ be the set
of all $x\in\bX$ such that the orbit closure $\overline{(T\times T)x}$
contains $y$. Then $\bX_y$ is an open affine $T\times T$-stable subset
of $\bX$, containing $y$ as its unique closed $T\times T$-orbit (it is
the image under $(1,w_0)$ of the affine chart $\cV$ defined in
\cite[\S 2]{S}.) Because $y\in\bX(w)$, the intersection
$\bX(w)\cap\bX_y:=\bX(w)_y$ is a non-empty open affine 
$T\times T$-stable subset of $\bX(w)$, containing $y$ as its unique
closed $T\times T$-orbit. It follows that the Picard group of
$\bX(w)_y$ is trivial (because $\bX(w)$ is normal), and also that any
regular invertible function on $\bX(w)_y$ is constant. Therefore, any
Cartier divisor on $\bX(w)$ is rationally equivalent to a unique
Cartier divisor with support in the complement
$\bX(w)\setminus\bX(w)_y$. Now the abelian group of Weil divisors with
support in $\bX(w)\setminus\bX(w)_y$ is free of finite rank. 

We now prove that any Cartier divisor $D$ on $\bX(w)$ which is
numerically equivalent to zero (that is, $(D\cdot C)=0$ for each
closed curve $C$ in $\bX(w)$) is rationally equivalent to
zero. Indeed, by \cite[19.3.3]{F}, there exists a positive integer $m$
such that $mD$ is algebraically equivalent to zero. But algebraic and
rational equivalence coincide for Cartier divisors on $\bX(w)$, by 
freeness of ${\rm Pic}\;\bX(w)$ and \cite[19.1.2]{F}. Thus, the class
of $mD$ in ${\rm Pic}\;\bX(w)$ is zero, and we conclude by freeness of
${\rm Pic}\;\bX(w)$ again.

For a line bundle $L$ on $\bX(w)$, define a weight $\lambda$ by
$$
\lambda:=(L\cdot C(\alpha_1))\omega_1+\cdots+
(L\cdot C(\alpha_r))\omega_r,
$$
so that $(L\cdot C(\alpha_i))=(\cL_{\bX}(\lambda)\cdot C(\alpha_i))$
for $1\leq i\leq r$. By Lemma \ref{curves}, it follows that
$(L\cdot C)=(\cL_{\bX}(\lambda)\cdot C)$ for all closed curves $C$
in $\bX(w)$. By the previous step, $L$ is isomorphic to
$res_{\bX(w)}\cL_{\bX}(\lambda)$. This proves 
that ${\rm Pic}\;\bX \cong {\rm Pic}\;\bX(w)$. 

For the remaining assertions of Theorem \ref{pic},
let $L=\cL_{\bX}(\lambda)$ be a line bundle on $\bX$ such
that $res_{\bX(w)}(L)$ is generated by its global sections
(resp. ample). Then $(L\cdot C)\geq 0$ (resp. $>0$) for any closed
curve $C$ in $\bX(w)$. 
Applying this to the curves $C(\alpha_i)$, one obtains that
$\lambda$ is dominant (resp. dominant regular). 

Conversely, for dominant $\lambda$, the line bundle 
$\cL_{\bX}(\lambda)$ admits a global section $\sigma_{\lambda}$ which
does not vanish identically on $Y$, by Lemma \ref{div}. Thus, the 
$\tG\times\tG$-translates of $\sigma_{\lambda}$ generate
$\cL_{\bX}(\lambda)$.

If moreover $\lambda$ is regular, then $\cL_{\bX}(\lambda)$ is ample
by \cite[\S 2]{S}.
\end{proof}

\medskip 
For any weight $\lambda$, the space
$$
H^0(\bX(w),\cL_{\bX(w)}(\lambda)):=H^0(\bX(w),\lambda)
$$ 
is a finite-dimensional $\tB\times \tB$-module. Its
$\tB\times\tB$-submodules
$$
H^0(\bX(w),\cL_{\bX}(\lambda)\otimes\cI_{Z(w)}^n)
:=F_nH^0(\bX(w),\lambda)
$$
(where $n\in\mN$) form a decreasing filtration, which we call the
{\it canonical filtration}. 
Since $H^0(\bX(w),\lambda)$ is finite dimensional and since 
$\bigcap_{n\geq 0} F_nH^0(\bX(w),\lambda) = 0$, this filtration is 
finite, that is, there exists an integer $n_0(\lambda)$ such that 
$F_nH^0(\bX(w),\lambda) = 0$ for $n > n_0(\lambda)$.

Let $\bn=(n_1,\ldots,n_r)\in\mN^r$ 
and let $|\bn| = n_1+\cdots+n_r$. Then
multiplication by the section $\sigma_1^{n_1}\cdots \sigma_r^{n_r}$
defines a map 
$$
\sigma^{\bn}:H^0(\bX(w),\lambda-n_1\alpha_1-\cdots-n_r\alpha_r)\to
H^0(\bX(w),\lambda).
$$
Because each $\sigma_i$ is $\tG\times \tG$-invariant and non identically
zero on $\bX(w)$, this map is injective and 
$\tB\times \tB$-equivariant. Let $F_{\bn}H^0(\bX(w),\lambda)$ be 
the image of $\sigma^{\bn}$; it is a $\tB\times\tB$-submodule of 
$F_n H^0(\bX(w),\lambda)$, where $n = |\bn|$.

\begin{theorem}\label{fil}
With notation as above, we have
$$
F_n H^0(\bX(w),\lambda)=
\sum_{ \deuxind{\bn\in\mN^r}{|\bn| = n} }
F_{\bn}H^0(\bX(w),\lambda)
$$

Further, the $n$-th layer of the associated graded module satisfies
$$
gr_n H^0(\bX(w),\lambda)=\bigoplus_{(n_1,\ldots,n_r)}
H^0(Z(w),\lambda-n_1\alpha_1-\cdots-n_r\alpha_r),
$$
the sum being taken
over all $(n_1,\ldots,n_r)\in\mN^r$ such that $n_1+\cdots+n_r=n$
and that $\lambda-n_1\alpha_1-\cdots-n_r\alpha_r$ is dominant.

In particular, 
$$
gr H^0(\bX(w),\lambda)=
\bigoplus_{ \deuxind{\mu\in\tcX^+}{\mu\leq\lambda} } 
H^0(Z(w),\mu). 
$$
\end{theorem}

\begin{proof}  
{F}rom the exact sequence of sheaves on $\bX(w)$:
$$
0\to \cI_{Z(w)}^{n+1}\otimes\cL_{\bX(w)}(\lambda)
\to\cI_{Z(w)}^n\otimes\cL_{\bX(w)}(\lambda)
\to\cI_{Z(w)}^n/\cI_{Z(w)}^{n+1}\otimes\cL_{\bX(w)}(\lambda)
\to 0,
$$
we see that $gr_n H^0(\bX(w),\lambda)$ injects into
$H^0(Z(w),
\cI_{Z(w)}^n/\cI_{Z(w)}^{n+1}\otimes\cL_{\bX(w)}(\lambda))$. 
The latter is equal to
$$
\bigoplus_{n_1+\cdots+n_r=n}
H^0(Z(w),\lambda-n_1\alpha_1-\cdots-n_r\alpha_r).
$$
We have indeed
$$
\cI_{Z(w)}^n/\cI_{Z(w)}^{n+1} = \bigoplus_{n_1+\cdots+n_r=n} 
\cL_{Z(w)}(-n_1\alpha_1-\ldots-n_r\alpha_r) 
\cdot \sigma_1^{n_1}\cdots \sigma_r^{n_r}
$$
because $\cI_{Z(w)}$ is generated by the regular sequence
$(\sigma_1,\ldots,\sigma_r)$ by Corollary \ref{reg}. We now need the
following 

\begin{lemma}\label{dom}
For a weight $\mu$, the following conditions are equivalent:

\noindent
(i) $\mu$ is dominant.

\noindent
(ii) $H^0(Z(w),\mu)$ is nonzero.
\end{lemma}

\begin{proof}
(i)$\Rightarrow$(ii) If $\mu$ is dominant, then the restriction to $Y$
of $\cL_{\bX}(\mu)$ is generated by its global sections.

(ii)$\Rightarrow$(i) Recall that
$$
Z(w)=\bigcup_{ \deuxind{x\in W}{\ell(wx)=\ell(w)+\ell(x)} } 
S(wx)\times S(xw_0), 
$$
and that the restriction of $\cL_{\bX}(\mu)$ to $Y$ is equal to
$\cL_Y(\mu,-w_0\mu)$. Thus, there exists $w\in W$ such that both
$H^0(S(wx),\mu)$ and $H^0(S(xw_0),-w_0\mu)$ are non-zero. But
$H^0(S(wx),\mu)\neq 0$ implies that
$\langle\mu,\check\alpha\rangle\geq 0$ for each $\alpha\in\Delta$ such
that $wx\alpha\in\Phi^-$; 
this follows from \cite[Cor. 2.3]{P}, see also \cite{Da}. 
Similarly, $H^0(S(xw_0),-w_0\mu)\neq 0$ implies that 
$\langle -w_0\mu,\check\beta\rangle\geq 0$ for each $\beta\in\Delta$
such that $xw_0\beta\in\Phi^-$. 
Since $-w_0$ permutes the simple roots, the latter is equivalent to 
$\langle\mu,\check\alpha\rangle\geq 0$ for each $\alpha\in\Delta$ such 
that $x\alpha\in\Phi^+$.  
Now, for each $\alpha\in\Delta$, we have either 
$x\alpha\in\Phi^+$ or $wx\alpha\in\Phi^-$, because
$\ell(wx)=\ell(w)+\ell(x)$. 
\end{proof}

Returning to the proof of Theorem \ref{fil},
let $\mu=\lambda-n_1\alpha_1-\cdots-n_r\alpha_r$ such that the space
$H^0(Z(w),\mu)$ is nonzero. Then $\mu$ is dominant by
Lemma \ref{dom}. By Corollary \ref{sur}, the restriction
$$
H^0(\bX(w),\mu)\to H^0(Z(w),\mu)
$$
is surjective; therefore, the restriction
$$
F_{\bn} H^0(\bX(w),\lambda)\to
H^0(Z(w),\lambda-n_1\alpha_1-\cdots-n_r\alpha_r)
$$
is surjective. It follows that, firstly, 
$$
gr_n H^0(\bX(w),\lambda) \cong \bigoplus_{(n_1,\ldots,n_r)}
H^0(Z(w),\lambda-n_1\alpha_1-\cdots-n_r\alpha_r),
$$
where the sum is taken 
over all $(n_1,\ldots,n_r)\in\mN^r$ such that $n_1+\cdots+n_r=n$
and that $\lambda-n_1\alpha_1-\cdots-n_r\alpha_r$ is dominant, 
and, secondly, that $F_{\bn} H^0(\bX(w),\lambda)$ is equal to
$$
\sum_{n_1+\cdots+n_r=n}F_{(n_1,\ldots,n_r)}H^0(\bX(w),\lambda)
+ F_{n+1}H^0(\bX(w),\lambda).
$$
Since $F_{n+1}H^0(\bX(w),\lambda) = 0$ for $n\gg 0$, 
this implies our statements.
\end{proof}

In particular, one obtains the following 

\begin{corollary}\label{effective} 
$\cL_{\bX(w)}(\lambda)$ admits a nonzero global section if and only if 
$\lambda$ belongs to the monoid generated by $\Delta$ and $\tcX^+$. 
\end{corollary}

Consider now the space
$$
R(w):=\bigoplus_{\lambda\in\tcX} H^0(\bX(w),\lambda). 
$$
Then $R(w)$ is a ring, with a grading by $\tcX$. 
By Theorem \ref{pic}, $R(w)$ can be seen as the
multihomogeneous coordinate ring of $\bX(w)$. Observe that
$\sigma_1,\ldots,\sigma_r$ are homogeneous elements of $R(w)$ of
degrees $\alpha_1,\ldots,\alpha_r$. 

We define similarly
$$
A(w):=\bigoplus H^0(Z(w),\mu)
$$
(sum over all weights $\mu$, or, equivalently, over all dominant
weights by Lemma \ref{dom}). Then $A(w)$ is the multihomogeneous
coordinate ring over $Z(w)$, a union of Schubert varieties 
in $G/B \times G/B$.

\begin{corollary}\label{ring}
The ring $R(w)$ is generated by its subspaces $H^0(\bX(w),\omega_i)$
$(1\leq i\leq r)$, together with $\sigma_1,\ldots,\sigma_r$. The
latter form a regular sequence in $R(w)$, and the quotient
$R(w)/(\sigma_1,\ldots,\sigma_r)$ is isomorphic to $A(w)$. 
\end{corollary}

\begin{proof}
The canonical filtrations of the $H^0(\bX(w),\lambda)$ fit
together into a filtration $(F_n R(w))$ of $R(w)$. Theorem \ref{fil}
implies that
$$
F_n R(w) = (\sigma_1,\ldots,\sigma_r)^n
$$
(the $n$-th power of the ideal generated by
$\sigma_1,\ldots,\sigma_r$) and that
$$
F_n R(w)/F_{n+1} R(w) = 
\bigoplus_{ \deuxind{n_1+\ldots+n_r=n}{\mu\in\tcX^+} } 
\sigma_1^{n_1}\cdots\sigma_r^{n_r} H^0(Z(w),\mu)
$$
Thus, the associated graded ring is isomorphic to the polynomial
ring $A(w)[t_1,\ldots,t_r]$. By \cite[1.1.15]{BH} 
$\sigma_1,\ldots,\sigma_r$ form a regular
sequence in $R(w)$. Further, by \cite[II.14.15, II.14.21]{J}, the graded
ring $A(w)$ is a quotient of $\bigoplus_{\mu} H^0(Y,\mu)$, and the
latter ring is generated by its subspaces $H^0(Y,\omega_i)$. So $R(w)$
is generated by  $\sigma_1,\ldots,\sigma_r$ and the
$H^0(\bX(w),\omega_i)$ (which lift the $H^0(Z(w),\omega_i)$.)
\end{proof}

We will show in Section 7 that the rings $R(w)$ and $A(w)$ are
Cohen-Macaulay.

\section{The van der Kallen filtration}

In this section, we construct geometrically van der Kallen's
filtration of the $\tB\times\tB$-module $k[\tB]$ (the ring of regular
functions on $\tB$), see \cite[Th. 1.13]{vdK1}. For this, we first
obtain a coarser filtration whose layers are spaces of global sections
of line bundles over $\oB\cap Y=Z$. In particular, $k[\tB]$ admits a
Schubert filtration as defined in \cite[2.8]{P} (see also
\cite[6.3.4]{vdK2}).

For $\mu\in\tcX^+$, we set 
$$
H^0(Z,\mu):=M(\mu).
$$
Then $M(\mu)$ is a finite dimensional $\tB\times\tB$-module. By
Theorem \ref{fil}, each $\tB\times\tB$-module $H^0(\oB,\lambda)$
has a filtration with layers $M(\mu)$ where $\mu\in\tcX^+$ and
$\mu\leq\lambda$. 

We will need more notation on $\tB$-modules, taken from \cite{vdK1}. 
Let $\nu$ be a weight; then there exist a unique $w=w_{min}\in W$ 
and a unique dominant weight $\mu$ such that $\nu=w\mu$ and that
the length of $w$ is minimal. Set
$$
P(\nu):=H^0(S(w_{min}),\mu) \text{ and }
Q(\nu):=H^0(S(w_{min}),\cI_{\partial S(w_{min})}\otimes\cL_{G/B}(\mu))
$$
where $\partial S(w)$ denotes the boundary of $S(w)$, that is, the
complement of its open $B$-orbit $BwB/B$. Then both $P(\nu)$ and
$Q(\nu)$ are finite dimensional $\tB$-modules. 

\begin{theorem}\label{vdk}
(i) The $\tB\times\tB$-module $k[\tB]$ has a canonical increasing
filtration by finite dimensional submodules, with associated graded
$$
\bigoplus_{\mu\in\tcX^+} M(\mu).
$$

\noindent
(ii) For any $\mu\in\tcX^+$, the $\tB\times\tB$-module
$M(\mu)$ has a filtration with associated graded
$$
\bigoplus_{\nu\in W\mu} P(\nu)\boxtimes Q(-\nu).
$$
\end{theorem}

\begin{proof} 
Set $\Gamma = \tcX/\cX$. For $\gamma\in\Gamma$, let 
$k[\tB]_{\gamma}$ be the sum of $\tT$-weight spaces in $k[\tB]$ (for
the right $\tT$-action) over all weights in the coset $\gamma$. Then
each $k[\tB]_{\gamma}$ is a $\tB\times\tB$-submodule of $k[\tB]$, and
we have
$$
k[\tB] = \bigoplus_{\gamma\in\Gamma} k[\tB]_{\gamma}.
$$
Further, the $k[B]$-module $k[\tB]_{\gamma}$ is freely generated by
any $B$-eigenvector. 

Choose $\gamma\in\Gamma$ and a dominant weight $\lambda$ in the coset
$\gamma$. By Theorem \ref{fil}, $H^0(\oB,\lambda)$ contains a
submodule isomorphic to $M(\nu)$, for some 
$\nu\in\tcX^+ \cap \gamma$. Thus, $H^0(\oB,\lambda)$ contains a right
$\tB$-eigenvector $v_\nu$ of weight $\nu$, and one deduces that 
$$
H^0(B,\cL_{\oB}(\lambda)) \cong k[B]\, v_{\nu} = k[\tB]_{\gamma}.
$$

Now, let us filter $k[\tB]_{\gamma}$ by the order of pole along the
boundary of $\oB$. Specifically, consider the section 
$\sigma:=\sigma_1\cdots\sigma_r$ of $\cL_{\oB}(\beta)$, where
$\beta=\alpha_1+\cdots+\alpha_r$. Then $\sigma$ is invariant by
$\tB\times\tB$, and its zero set is the boundary $\oB-B$. Therefore,
the $\tB\times\tB$-module $k[\tB]_{\gamma}=H^0(B,\cL_{\oB}(\lambda))$
is the increasing union of its finite dimensional submodules
$$
\sigma^{-n} H^0(\oB,\lambda+n\beta)
$$
for $n\geq 0$. The associated graded of this filtration satisfies
$$
gr^n k[\tB]_{\lambda} \cong
H^0(\oB,\lambda+n\beta)/\sigma H^0(\oB,\lambda+(n-1)\beta).
$$
Let $R$ be the multihomogeneous coordinate ring on $\oB$, then
$\sigma\in R_{\beta}$ and 
$$
gr^n k[\tB]_{\lambda} \cong
R_{\lambda+n\beta}/\sigma R_{\lambda+(n-1)\beta}=
(R/\sigma R)_{\lambda+n\beta}.
$$
Consider the decreasing filtration of $R/\sigma R$, image of 
the filtration of $R$ by the ideals 
$(\sigma_1,\ldots,\sigma_r)^m R$. As $\sigma_1,\ldots,\sigma_r$ 
form a homogeneous regular sequence in $R$, the associated graded of 
$R/\sigma R$ satisfies 
$$
gr_m R/\sigma R=
\bigoplus_{\deuxind{m_1+\cdots+m_r=m}{m_1\cdots m_r=0}} 
\sigma_1^{m_1}\cdots \sigma_r^{m_r} A
$$
where $A$ is the multihomogeneous coordinate ring of $Z$. Taking
homogeneous components of degree $\lambda+n\beta$, we see that
each $gr^n k[\tB]_{\lambda}$ has a 
finite decreasing filtration with associated graded
$$
\bigoplus_{\deuxind{(m_1,\ldots,m_r)\in\mN^r}
{\lambda+n\beta-m_1\alpha_1-m_r\alpha_r\in\tcX^+}}
M(\lambda+n\beta-m_1\alpha_1-m_r\alpha_r).
$$
Reordering the indices, we obtain a canonical increasing filtration of
$k[\tB\times\tB]$ satisfying the requirements of (i).

For (ii), recall that the irreducible components of $Z$ are exactly
the $S(w)\times S(ww_0)$ for $w\in W$. We first construct an
increasing filtration of $Z$ by partial unions of these components, as
follows. Choose an indexing 
$W=\{w_1,\ldots,w_M\}$ which is compatible with the Bruhat-Chevalley
order, that is, $i\leq j$ if $w_i\leq w_j$. In particular, $w_1=1$ and
$w_M=w_0$. Set  
$$
Z_i:=S(w_i)\times S(w_iw_0),~Z_{\geq i}:=\bigcup_{j\geq i}Z_j,~
Z_{>i}=\bigcup_{j>i} Z_j.
$$
Then we have the following

\begin{lemma}\label{int}
$Z_i\cap Z_{>i}=S(w_i)\times\partial S(w_iw_0)$.
\end{lemma}

\begin{proof}
Let $x$, $y$ in $W$ such that $S(x)\times S(yw_0)\subseteq Z_i$, that
is, $x\leq w_i\leq y$. If moreover 
$S(x)\times S(yw_0)\subseteq Z_{>i}$, then $w_j\leq y$ for some
$j>i$. It follows that $y\neq w_i$, so that 
$Z_i\cap Z_{>i}\subseteq S(w_i)\times \partial S(w_iw_0).$

For the opposite inclusion, let $y\in W$ such that
$S(yw_0)\subset\partial S(w_iw_0)$, that is, $y>w_i$. Then $y=w_j$ for
some $j>i$. Thus, $S(w_i)\times S(yw_0)\subseteq S(y)\times S(yw_0)$
is contained in $Z_i\cap Z_{>i}$.
\end{proof}

Returning to the proof of Theorem \ref{vdk}, let $\cI_i$ be the ideal
sheaf of $Z_{>i}$ in $Z_{\geq i}$. Then $\cI_i$ identifies to the
ideal sheaf of $Z_i\cap Z_{>i}$ in $Z_i$. By definition, we have an exact
sequence of sheaves of $\cO_Z$-modules:
$0\to \cI_i \to \cO_{Z_{\geq i}} \to \cO_{Z_{>i}} \to 0$. Thus, the
sequence 
$$
0\to\cI_i\otimes\cL_Z(\mu)\to\cL_{Z_{\geq i}}(\mu)\to\cL_{Z_{>i}}(\mu)
\to 0
$$ 
is exact. Further, $H^1(Z_{\geq i},\cL_Z(\mu))=0$ as $Z_{\geq i}$ is
a union of Schubert varieties in $Y$, and $\mu$ is dominant. 
So we obtain an exact sequence
$$
0\to H^0(Z_{\geq i},\cI_i\otimes\cL_Z(\mu)) \to
H^0(Z_{\geq i},\mu) \to H^0(Z_{> i},\mu) \to 0.
$$
Now Lemma \ref{int} implies that
$$\displaylines{
H^0(Z_{\geq i},\cI_i\otimes\cL_Z(\mu))
=H^0(Z_i,\cI_i\otimes\cL_Z(\mu))
\hfill\cr\hfill
=H^0(S(w_i),\mu)\boxtimes 
H^0(S(w_iw_0),\cI_{\partial S(w_iw_0)}\otimes\cL_{G/B}(-w_0\mu)).
\cr}$$
By induction on $i$, we thus obtain a filtration of
$M(\mu)$ with associated graded
$$\bigoplus_{x\in W}
H^0(S(x),\mu)\boxtimes 
H^0(S(w_0x), \cI_{\partial S(w_0x)}\otimes\cL_{G/B}(-w_0\mu)).
$$
Further, we have $H^0(S(x),\mu)=P(x\mu)$ by 
\cite[Lemma 2.3.2]{vdK2}. And  
$$
H^0(S(x),\cI_{\partial S(x)}\otimes\cL_{G/B}(\mu))=Q(x\mu)
$$ 
if $x$ is the element of minimal length in its coset $xW_{\mu}$ (where
$W_{\mu}$ is the isotropy group of $\mu$ in $W$.) Otherwise, we claim
that $H^0(S(x),\cI_{\partial S(x)}\otimes\cL_{G/B}(\mu))=0$. Indeed,
by [{\sl loc. cit.}], the restriction map
$H^0(S(x),\mu)\to H^0(S(x_{min}),\mu)$
is an isomorphism.

It follows that the $\tB$-module
$H^0(S(xw_0),\cI_{\partial S(xw_0)}\otimes\cL_{G/B}(-w_0\mu))$ equals
$Q(-x\mu)$ if $x$ has maximal length in its coset $xW_{\mu}$, and
equals $0$ otherwise. Thus, the $\tB\times\tB$-module
$M(\mu)$ has a filtration with associated graded
$\oplus P(x\mu)\boxtimes Q(-x\mu)$, sum over all $x\in W$ such that
$x$ has maximal length in its $W_{\mu}$-coset. This implies (ii). 
\end{proof}

\smallskip\noindent
{\sc Remark.} A similar argument shows that the 
$\tG\times\tG$-module $k[\tG]$ has an increasing filtration by finite
dimensional submodules with associated graded 
$$
\bigoplus_{\mu\in\tcX^+} H^0(G/B,\mu)\boxtimes H^0(G/B,-w_0\mu).
$$
This gives a geometric proof of a result of Donkin and Koppinen
\cite[II.4.20]{J}. 

\medskip 
For any dominant weight $\mu$, we denote by $c_{\mu}={\rm ch}\;M(\mu)$
the character of the finite-dimensional $\tT\times\tT$-module
$M(\mu)$. Then $c_{\mu}$ is a regular function on
$\tT\times\tT$, and we have by Theorem \ref{fil}, for $t,u\in \tT$:
$$
c_{\mu}(t,u)=\sum_{\nu\in W\mu} ch\;P(\nu)(t)\; ch\;Q(-\nu)(u). 
$$
Further, $ch\;P(\nu)$ is given by the Demazure character formula
\cite[II.14.18]{J}, and $ch\;Q(\nu)$ is given by a closely related
formula \cite[Th. 2.1]{M}.

By Corollary \ref{sur}, we have $H^i(Z,\mu)=0$ for
$i\geq 1$; thus, we can extend the map
$\mu\mapsto c_{\mu}$ to the group $\tcX$ by setting 
$$
c_{\lambda}=\chi(Z,\lambda)
=\sum_{i\geq 0} (-1)^i {\rm ch}\;H^i(Z,\lambda)
$$
for arbitrary $\lambda\in\tcX$. 

We now establish two symmetry properties of the resulting map
$\lambda\mapsto c_{\lambda}$; the second symmetry will be an essential
ingredient of the proof that $\oB$ is Gorenstein. We will determine
the value of $c_{\lambda}$ at $(t,t^{-1})$ and, in particular, the
dimension of $M(\lambda)$ in Corollary \ref{sep} below.

\begin{theorem}\label{recip}
We have
$c_{-w_0\lambda}(t,u)=c_{\lambda}(u,t)$ and
$$
c_{-\lambda}(t^{-1},u^{-1})=(-1)^N\rho(t)\rho(u)
c_{\lambda-\rho}(t,u)
$$
for all $\lambda\in\tcX$ and $t,u\in\tT$.
\end{theorem}

\begin{proof} 
With notation as in Lemma \ref{int}, we have for any
$\lambda\in\tcX$ and any index $i$:
$$\displaylines{
\chi(Z_{\geq i},\lambda)=
\chi(Z_{>i},\lambda)+\chi(Z_{\geq i},\cI_i\otimes\cL_Z(\lambda))
\hfill\cr\hfill
=\chi(Z_{>i},\lambda)+\chi(S(w_i),\lambda)\, 
\chi(S(w_iw_0),\cI_{\partial S(w_iw_0)}\otimes\cL_{G/B}(-w_0\lambda))
\cr}$$
by Lemma \ref{int} and the argument thereafter. 
Since $Z = Z_{\geq 1}$, it follows that
$$
\chi(Z,\lambda)=\sum_{x\in W}\chi(S(x),\lambda)\, 
\chi(S(xw_0),\cI_{\partial S(xw_0)}\otimes\cL_{G/B}(-w_0\lambda)) \eqno (*)
$$
for all $\lambda\in\tcX$. 

Recall now that each Schubert variety $S(x)$ is Cohen-Macaulay. Denote
by $\omega_{S(x)}$ its canonical sheaf, a $B$-linearized
sheaf. By \cite[A.4]{vdK2}, we have an isomorphism of $B$-linearized
sheaves
$$
\omega_{S(x)}\cong\cI_{\partial S(xw_0)}\otimes\cL_{G/B}(-\rho)[\rho]
$$
where $[\rho]$ denotes the shift by the character $\rho$ in the
$B$-linearization. Thus, we obtain by using Serre duality on each
$S(x)$ to pass from the first to the second line: 
$$\displaylines{
c_{-\lambda}(t^{-1},u^{-1})=\sum_{x\in W}
\chi(S(x),-\lambda)(t^{-1})\;
\chi(S(xw_0),\cI_{\partial S(xw_0)}\otimes\cL_{G/B}(w_0\lambda))(u^{-1})
\cr
=(-1)^N\rho(t)\rho(u)\sum_{x\in W}
\chi(S(x),\cI_{\partial S(x)}\otimes\cL_{G/B}(-\rho+\lambda))(t)\; 
\chi(S(xw_0),-\rho-w_0\lambda)(u)
\cr
=(-1)^N\rho(t)\rho(u)\sum_{y\in W}
\chi(S(y),-\rho-w_0\lambda)(u)\,
\chi(S(yw_0),\cI_{\partial S(yw_0)}\otimes\cL_{G/B}(-\rho+\lambda))(t)
\cr
=(-1)^N\rho(t)\rho(u)c_{-w_0\lambda-\rho}(u,t).
\cr}$$

On the other hand, set $Z^i:=S(w_i)\times S(w_iw_0)$ and define
similarly $Z^{\geq i}$, $Z^{>i}$. Then we obtain as in Lemma \ref{int} that  
$Z^i\cap Z^{>i}=\partial S(w_iw_0)\times S(w_i)$. As above, it follows
that
$$
\chi(Z,\lambda)=\sum_{x\in W}
\chi(S(xw_0),\cI_{\partial S(xw_0)}\otimes\cL_{G/B}(\lambda))\;
\chi(S(x),-w_0\lambda)
$$ 
for any $\lambda$. Thus,
$\chi(Z,-w_0\lambda)(u,t)=\chi(Z,\lambda)(t,u)$ and the first identity
is proved. In particular,
$c_{-w_0\lambda-\rho}(u,t)=c_{\lambda-\rho}(t,u)$ 
which completes the proof of the second identity.
\end{proof} 

\section{Closures of Borel subgroups are Gorenstein}

Let $R$ (resp. $A$) be the multihomogeneous coordinate ring of $\oB$
(resp. $Z=\oB\cap Y$) as defined in Section 2. We show that both $R$
and $A$ are Gorenstein; as a consequence, $\oB$ and $Z$ are Gorenstein
as well. 

It will be convenient to set
$$
\beta:=\alpha_1+\cdots+\alpha_r.
$$
Then the canonical sheaf of $\bX$ is $\cL_{\bX}(-2\rho-\beta)$ 
by \cite[\S3]{S}.

\begin{theorem}\label{gor} (i) $\oB$ $($resp. $Z)$ is
Gorenstein with canonical sheaf $\cL_{\oB}(-\rho-\beta)[\rho,\rho]$
$($resp. $\cL_Z(-\rho)[\rho,\rho])$ as a $B\times B$-linearized 
sheaf.

\noindent
(ii) The graded ring $R$ $($resp. $A)$ is Gorenstein and its canonical
module is generated by a homogeneous element of degree $\rho+\beta$
$($resp. $\rho)$, eigenvector of $B\times B$ of weight $(\rho,\rho)$.
\end{theorem}

\begin{proof} 
We begin by proving that $Z$ is Cohen-Macaulay. For this, we use the
notation introduced in the proof of Theorem \ref{vdk}. We check by
decreasing induction on $i$ that each $Z_{\geq   i}$ is 
Cohen-Macaulay. If $i=M$ then $Z_{\geq M}=S(w_0)\times S(1)\cong G/B$,
a non-singular variety. For arbitrary $i$, we have an exact sequence
$$
0\to\cO_{Z_{\geq i}}\to\cO_{Z_i}\oplus \cO_{Z_{> i}}\to
\cO_{Z_i\cap Z_{>i}}\to 0.
$$
Further, we know that $Z_i=S(w_i)\times S(w_iw_0)$ is Cohen-Macaulay;
and, by the induction hypothesis, the same holds for $Z_{>i}$. On the
other hand, $Z_i\cap Z_{>i}=S(w_i)\times \partial S(w_iw_0)$ by Lemma
\ref{int}. The canonical sheaf of $S(w_iw_0)$ is the tensor
product of the ideal sheaf of $\partial S(w_iw_0)$ with the invertible
sheaf $\cL_{G/B}(-\rho)$. By \cite[Proposition 3.3.18]{BH}, it follows
that $\partial S(w_iw_0)$ is Cohen-Macaulay, of depth 
$\ell(w_iw_0)-1$. Thus, the depth of $Z_i\cap Z_{>i}$ is
$\ell(w_i)+\ell(w_iw_0)-1=\ell(w_0)-1=\dim Z_{\geq i} -1$. 
Together with the exact sequence above, this implies easily that
$Z_{\geq i}$ is Cohen-Macaulay, see \cite[Proposition 1.2.9]{BH}.

We now prove that the ring $A$ is Cohen-Macaulay. For this, let
$C={\rm Spec}\;A$ be the corresponding affine scheme. Then $C$ is the
multicone over $Z$ in the sense of \cite{KR}; we now recall some
constructions from that paper. Let $E$ be the total space of the
vector bundle over $Z$, equal to the direct sum of the line bundles
$\cL_Z(-\omega_1),\ldots,\cL_Z(-\omega_r)$; let $q:E\to Z$ be 
the projection. Then 
$$
q_*\cO_E=Sym_{\cO_Z}\bigoplus_{i=1}^r\cL_Z(\omega_i)
=\bigoplus_{\mu\in\tcX^+}\cL_Z(\mu).
$$
In particular, $H^0(E,\cO_E)=A$
so that we have a morphism $p:E\to C$. The torus $\tT$ acts on
$E$ and on $C$ (because $A$ is graded by the character group of
$\tT$), compatibly with the action of $\tB\times\tB$. 
Clearly, $p$ and $q$ are equivariant for the action of
$\tB\times\tB\times\tT$.

As the line bundles $\cL_Z(\omega_1),\ldots,\cL_Z(\omega_r)$ are
generated by their global sections, $p$ is proper. Further, we have
$p_*\cO_E=\cO_C$ as $C$ is affine and $H^0(E,\cO_E)=H^0(C,\cO_C)$. In
particular, $p$ is surjective. 

Let $E^0$ be the total space of $E$ minus the union of all sub-bundles
$\oplus_{j\neq i}\cL_Z(\omega_j)$ for $1\leq i\leq r$; let 
$p^0:E^0\to C$ and $q^0:E^0\to Z$ be the restrictions of $p$ and
$q$. Then $p^0$ is an isomorphism onto an open subset $C^0$ of $C$, and
$q^0$ is a principal $\tT$-bundle over $Z$. As a consequence,
the restriction of $p$ to each irreducible component of $E$ (that is,
to each  $q^{-1}(S(w)\times S(ww_0))$ is birational. Thus,
$C$ is equidimensional of dimension $\dim(B)=N+r$.

We claim that $R^ip_*\cO_E=0$ for $i\geq 1$. Indeed, as 
$C$ and $q$ are affine, this amounts to:
$$
0=H^i(E,\cO_E)=H^i(Z,q_*\cO_E)=\bigoplus_{\mu\in\tcX^+} H^i(Z,\mu), 
$$
which follows from Corollary \ref{sur}.

Because $Z$ is Cohen-Macaulay, the same holds for $E$, and we have
$$
\omega_E=q^*\omega_Z\otimes\omega_{E/Z}=q^*(\omega_Z\otimes
\cL_Z(\omega_1)[-\omega_1]\otimes\cdots
\otimes\cL_Z(\omega_r)[-\omega_r])
=q^*(\omega_Z\otimes\cL_Z(\rho))[-\rho]
$$
as a $\tB\times\tB\times\tT$-linearized sheaf; here $[\chi]$ denotes
the shift of the $\tT$-linearization by the character $\chi$.
We claim that $R^ip_*\omega_E=0$ for $i\geq 1$, that is,
$H^i(E,\omega_E)=0$ for $i\geq 1$. Indeed, we have
$$
H^i(E,\omega_E)= H^i(Z,q_*\omega_E)
=\bigoplus_{\mu\in\tcX^+} H^i(Z,\omega_Z\otimes\cL_Z(\rho+\mu)) 
=\bigoplus_{\mu\in\tcX^+} H^{N-i}(Z,-\rho-\mu)^* ,
$$
as $Z$ is equidimensional of dimension $N$. 
For $\mu\in\tcX^+$, the line bundle $\cL_Z(\rho+\mu)$ is
ample. Because $Z$ is Cohen-Macaulay, we have therefore 
$$
H^j(Z,-n(\rho+\mu))=0
$$ 
for $j<N$ and large $n$. But $Z$, being an union of Schubert varieties
in $Y$, is Frobenius split. Thus, $H^j(Z,-\rho-\mu)=0$ for $j<N$ by
\cite[Proposition A.2.1]{vdK2}. This proves our claim.

We now recall a version of a result of Kempf, see e.g. 
\cite[p. 49-51]{KKMS}. 

\begin{lemma}\label{dual}
Let $p:\hat X\to X$ be a proper morphism of algebraic
schemes. Assume that $\hat X$ is Cohen-Macaulay, 
$X$ is equidimensional of the same dimension as $\hat X$,
$p_*\cO_{\hat X}=\cO_X$ and 
$R^ip_*\cO_{\hat X}=R^ip_*\omega_{\hat X}=0$ for $i\geq 1$. Then
$X$ is Cohen-Macaulay with dualizing sheaf $p_*\omega_{\hat X}$.
\end{lemma}

\begin{proof} 
The statement is local in $X$, so that we may assume that $X$ is a
closed subscheme of a smooth affine scheme $S$. Denote by 
$\iota:X\to S$ the inclusion and set $\pi:=\iota\circ p$. Then 
$\pi_*\cO_{\hat X}=\iota_*\cO_X$ and $R^i\pi_*\cO_{\hat X}=0$ for
$i\geq 1$. Applying the duality theorem to the proper morphism
$\pi:\hat X\to S$ and the sheaves $\cO_{\hat X}$ and $\omega_S$,
we obtain
$$\displaylines{
RHom_S(\iota_*\cO_X,\omega_S)=RHom_S(R\pi_*\cO_{\hat X},\omega_S)=
R\pi_*RHom(\cO_{\hat X},\pi^!\omega_S)
\hfill\cr\hfill
=R\pi_*\pi^!\omega_S=R\pi_*\omega_{\hat X}[\dim X-\dim S]
=\pi_*\omega_{\hat X}[\dim X-\dim S],
\cr}$$
that is, $Ext^i_S(\iota_*\cO_X,\omega_S)=0$ for $i\neq\dim S-\dim X$,
and 
$$
Ext^{\dim X-\dim S}_S(\iota_*\cO_X,\omega_S)=\pi_*\omega_{\hat X}.
$$
This means that $X$ is Cohen-Macaulay with canonical sheaf
$p_*\omega_{\hat X}$.
\end{proof}

Lemma \ref{dual} implies that the graded ring $A$ is Cohen-Macaulay
with canonical module
$$
\omega_A=H^0(C,p_*\omega_E)
=\bigoplus_{\mu\in\tcX^+} H^0(Z,\omega_Z\otimes\cL_Z(\mu+\rho))[-\rho]
=\bigoplus_{\mu\in\tcX^+} H^N(Z,-\mu-\rho)^*[-\rho].
$$
Further, $H^j(Z,-\mu-\rho)=0$ for $j\neq N$.
The module $\omega_A$ is $\tcX$-graded and each homogeneous 
component is a finite-dimensional $\tB\times\tB$-module. Thus, we can
consider the Hilbert series
$$
H_{\omega_A}(t,u,z):=\sum_{\lambda\in\tcX} 
ch\;\omega_{A,\lambda}(t,u)\;z^{\lambda}
$$ 
where $t$, $u$ are in $\tT$, and the $z^{\lambda}$ are the canonical
basis of the group algebra $\mZ[\tcX]$. Now we have
$$\displaylines{
H_{\omega_A}(t,u,z)=\sum_{\mu\in\tcX^+}
ch H^N(Z,-\mu-\rho)(t^{-1},u^{-1})\;z^{\mu+\rho}
\hfill\cr\hfill
=(-1)^N\sum_{\mu\in\tcX^+}
c_{-\mu-\rho}(t^{-1},u^{-1})z^{\mu+\rho}.
\cr}$$
Together with Theorem \ref{recip}, it follows that
$$
H_{\omega_A}(t,u,z) = \rho(t) \rho(u) z^{\rho} H_A(t,u,z).
$$ 
Using \cite[Cor. 4.3.8.a)]{BH}, we have therefore
$$
H_A(t^{-1},u^{-1},z^{-1}) = 
(-1)^{\dim(A)} \rho(t)\rho(u) z^{\rho} H_A(t,u,z).
$$

Now a result of Stanley \cite[Cor. 4.3.8.c)]{BH} would imply that
$A$ is Gorenstein if $A$ were a domain. This is not the case, but $A$
is the quotient of the domain $R$ by the ideal generated by the
regular sequence $(\sigma_1,\ldots,\sigma_r)$. It follows that $R$ is
Cohen-Macaulay with Hilbert series
$$
H_R(t,u,z)=\frac{H_A(t,u,z)}{\prod_{i=1}^r(1-z^{\alpha_i})}
$$
because each $\sigma_i$ is the restriction of a
$\tG\times\tG$-invariant section. As a consequence, we obtain
$$
H_R(t^{-1},u^{-1},z^{-1})
=(-1)^{\dim(R)}\rho(t)\rho(u)z^{\rho+\beta}H_R(t,u,z).
$$
Thus, by the result of Stanley quoted above, $R$ is Gorenstein and its 
canonical module is generated by a homogeneous element of degree
$\rho+\beta$, eigenvector of $\tB\times\tB$ of weight
$(\rho,\rho)$. It follows that $A$ is Gorenstein as well and that its
canonical module is generated in degree $\rho$ and weight $(\rho,\rho)$.

It remains to prove that $\oB$ and $Z$ are Gorenstein and to determine
their canonical sheaves. For this, consider the isomorphism
$p^0:E^0\to C^0$ where $C^0$ is an open subset of $C$, and the
principal $\tT$-bundle $q^0:E^0\to Z$. Then
$\omega_{E^0}=\cO_{E^0}[\rho,\rho]$  
as a $\tB\times\tB$-linearized sheaf, because the same holds for
$C$. Further, $\omega_{E^0}=q^{0*}(\omega_Z\otimes\cL_Z(\rho))$ 
so that
$$
\omega_Z\otimes\cL_Z(\rho)\otimes q_{0*}\cO_{E^0}
=q_{0*}\cO_{E^0}[\rho,\rho].
$$
Taking invariants of $\tT$, we obtain 
$$
\omega_Z\otimes\cL_Z(\rho)=\cO_Z[\rho,\rho],
$$
that is, $Z$ is Gorenstein with canonical sheaf
$\cL_Z(-\rho)[\rho,\rho]$. The argument for $\oB$ is similar.
\end{proof}

\section{The class of the diagonal for flag varieties}

We will construct a degeneration of the diagonal in $G/B\times G/B$
into a union of Schubert varieties. For this, we recall a special case
of a construction of \cite[1.6]{B}.

Consider the action of $B\times B$ on $\oB$ and the associated fiber
bundle 
$$
p:G\times G\times_{B\times B}\oB\to G/B\times G/B,
$$ 
a locally trivial fibration with fiber $\oB$. The action map 
$$G\times G\times\oB\to\bX=(G\times G)\oB:~(g,h,x)\mapsto (g,h)x$$
defines a $G\times G$-equivariant map
$$
\pi: G\times G\times_{B\times B}\oB\to \bX.
$$
Observe that $\pi$ factors through the closed embedding
$$
G\times G\times_{B\times B}\oB\to G/B\times G/B\times\bX:
(g,h,x)(B\times B)\mapsto(gB,hB,(g,h)x)
$$
followed by the projection
$G/B\times G/B\times X\to G/B\times G/B$. Thus, $\pi$ is proper and
its scheme-theoretic fibers identify to closed subschemes of 
$G/B\times G/B$ via $p_*$. Further, the fiber $\pi^{-1}(1)$ at the
identity is the diagonal $diag(G/B)$; and the reduced fiber
$\pi^{-1}(y)_{red}$ at the base point $y$ of $Y=G/B\times G/B$ is
$$
\bigcup_{x\in W} S(x)\times S(w_0x).
$$

Consider now the closure $\oT$ of $T$ in $\bX$, then $\oT$ is a
$T\times T$-stable subvariety fixed pointwise by $diag(T)$. By
\cite{S}, $\oT$ is smooth and meets $Y$ transversally at the points
$(w,ww_0)y$ for $w\in W$. Each of these points admits a 
$T\times T$-stable neighborhood isomorphic to affine $r$-space where
$T\times T$ acts linearly. Therefore, we can find a smooth curve
$\gamma\subseteq\oT$ isomorphic to affine line, containing $1$ 
(the identity element of $G$) and transversal to $Y$ at
$z:=(1,w_0)y$. In particular, $\gamma\setminus\{z\}$ is contained in
$T$. Further, $\pi^{-1}(\gamma)$ is a $diag(T)$-stable subvariety of 
$G\times G\times_{B\times B}\oB$ and we have a $diag(T)$-equivariant
isomorphism
$$
\pi^{-1}(\gamma)\simeq\{(gB,hB,x)\in G/B\times G/B\times\gamma
~\vert~ (g^{-1},h^{-1})x\in\oB\}
$$
identifying $\pi^{-1}(\gamma)\to\gamma$ to the restriction of the
projection $G/B\times G/B\times \gamma\to\gamma$.

\begin{theorem}\label{flat}
The morphism $\pi:G\times G\times_{B\times B}\oB\to\bX$ is flat, with
reduced fibers. Its restriction $\pi^{-1}(\gamma)\to\gamma$ is flat
and $diag(T)$-invariant, with fibers over $\gamma\setminus\{z\}$
isomorphic to $diag(G/B)$, and with fiber at $z$ equal to
$$
\bigcup_{x\in W} S(x)\times w_0S(w_0x).
$$
\end{theorem}

\begin{proof} 
By \cite[Proposition 1.6]{B}, $\pi$ is equidimensional. Further,
$G\times G\times_{B\times B}\oB$ is Cohen-Macaulay, as $\oB$
is. Because $\bX$ is smooth, it follows that $\pi$ is flat. 

For $x\in \bX$, the scheme-theoretic fiber $\pi^{-1}(x)$ identifies to 
$$
\{(gB,hB)\in G/B\times G/B~\vert~ (g^{-1},h^{-1})x\in\oB\}.
$$
Set $F:=\{(g,h)\in G\times G~\vert~ (g^{-1},h^{-1})x\in\oB\}$.
Then $F$ is stable under right multiplication by 
$B\times B$, and left multiplication by $(G\times G)_x$ (the isotropy
group of $x$ in $G\times G$). The quotient of $F$ by the right 
$B\times B$-action is $\pi^{-1}(x)$, whereas the quotient by the left
$(G\times G)_x$-action is isomorphic to the scheme-theoretic
intersection of $\oB$ with the orbit $(G\times G)\cdot x$. This
intersection is reduced by Corollary \ref{sur}; thus, $F$ and
$\pi^{-1}(x)$ are reduced, too.

The remaining asssertions are direct consequences of these facts.
\end{proof}

We now deduce from Theorem \ref{flat} a formula for the class in
equivariant $K$-theory of the diagonal of the flag variety.
Consider the diagonal action of $T$ on $G/B\times G/B$, and let
$K^T(G/B\times G/B)$ be the corresponding Grothendieck group of
$T$-linearized coherent sheaves. Then $K^T(G/B\times G/B)$ is a module
over the representation ring $R(T)$ of $T$; further, $W$ acts on
$K^T(G/B\times G/B)$ compatibly with its action on $R(T)$.
 
For a $T$-stable subvariety $S$ of $G/B\times G/B$, the class in
$K^T(G/B\times G/B)$ of the structure sheaf $\cO_S$ will be denoted by
$[S]$. In particular, we have the classes of Schubert varieties and of
their translates by $W\times W$; we also have the class of the
diagonal $diag(G/B)$. We will express the latter in terms of the
former.

This will imply a formula for the class of the diagonal in the
Grothendieck group $K(G/B\times G/B)$ of coherent sheaves on that
space, by applying the forgetful map
$$
K^T(G/B\times G/B)\to K(G/B\times G/B).
$$
Observe that the action of $G\times G$ on $K(G/B\times G/B)$
is trivial, because $G$ is generated by subgroups isomorphic to the
additive group.

To simplify our statements, we set for $x\in W$:
$$
S^-(x):=\overline{B^-wB}/B=w_0S(w_0x)
$$
and
$$
[S^-(x)]^0:=[S^-(x)]-[\partial S^-(x)]
=w_0[S(w_0x)]-w_0[\partial S(w_0x)].
$$

\begin{corollary}\label{diag}
With notation as above, we have in $K^T(G/B\times G/B)$:
$$\displaylines{
[diag(G/B)]=[\bigcup_{x\in W} S(x)\times S^-(x)]
\hfill\cr\hfill
=\sum_{x\in W} [S(x)]\times [S^-(x)]^0
=\sum_{x,y\in W, \,x\leq y} (-1)^{\ell(y)-\ell(x)}
[S(x)]\times [S^-(y)].
\cr}$$
As a consequence, we have in $K(G/B\times G/B)$:
$$\displaylines{
[diag(G/B)]=[\bigcup_{x\in W} S(x)\times S(w_0x)]
\hfill\cr\hfill
=\sum_{x\in W} [S(x)]\times [S(w_0x)]^0
=\sum_{x,y\in W,\,x\leq y} (-1)^{\ell(y)-\ell(x)}
[S(x)]\times[S(w_0y)].
\cr}$$
\end{corollary}

\begin{proof} 
As the map $\pi^{-1}(\gamma)\to\gamma$ is flat and
$diag(T)$-invariant, and $\gamma$ is isomorphic to affine line, the 
fibers $\pi^{-1}(1)$ and $\pi^{-1}(z)$ have the same class in
$K^T(G\times G\times_{B\times B}\oB)$. Thus, the direct images of
these fibers under $p$ are equal in $K^T(G/B\times G/B)$. 
This proves the first equality. 

For the second one, we use the notation introduced in the proof of
Theorem \ref{vdk}. Set 
$$
\Delta_i:=S(w_i)\times S^-(w_i),~
\Delta_{\geq i}:=\bigcup_{j\geq i}\;\Delta_i,~
\Delta_{> i}:=\bigcup_{j>i}\;\Delta_i.
$$
Then $\Delta_{\geq M}=S(w_0)\times w_0S(1)$, 
$\Delta_{\geq 1}=\bigcup_{x\in W}S(x)\times S^-(x)$,
and we obtain as in Lemma \ref{int} that:
$$
\Delta_i\cap\Delta_{>i}=S(w_i)\times \partial S^-(w_i).
$$
As a consequence, we have an exact sequence of sheaves
$$
0\to I_i \to \cO_{\Delta_{\geq i}}\to \cO_{\Delta_{>i}}\to 0
$$
where $I_i$ fits into an exact sequence
$$
0\to I_i\to\cO_{\Delta_i}\to\cO_{S(w_i)\times \partial S^-(w_i)}
\to 0.
$$
It follows that
$$\displaylines{
[\Delta_{\geq i}]=[\Delta_{>i}]+[I_i]
=[\Delta_{>i}]+[S(w_i)\times S^-(w_i)]
-[S(w_i)\times \partial S^-(w_i)]
\hfill\cr\hfill
=[\Delta_{>i}]+[S(w_i)]\times [S^-(w_i)]^0.
\cr}$$
By decreasing induction on $i$, we thus have 
$\Delta_{\geq i}=\sum_{j\geq i}\;[S(w_j)]\times [S^-(w_j)]^0$.

For the third equality, it suffices to prove that
$$
[S(x)]^0=\sum_{y\leq x} (-1)^{\ell(x)-\ell(y)} [S(y)]
$$
in $K^T(G/B)$. But the definition of $[S(x]^0$ implies that
$[S(x)]=\sum_{y\in W,\;y\leq x} [S(y)]^0$. Further, the M\"obius
function of the partially ordered set $(W,\leq)$ is given by
$\mu(y,x)=(-1)^{\ell(x)-\ell(y)}$ if $y\leq x$ and $\mu(y,x)=0$
otherwise, see \cite{D}.
\end{proof}

\smallskip 
Consider now the Grothendieck group $K^T(G/B)$ of $T$-linearized
coherent sheaves on $G/B$. Because $G/B$ is smooth and projective,
$K^T(G/B)$ is isomorphic to the Grothendieck group of $T$-linearized
vector bundles; as a consequence, it has the structure of a
$R(T)$-algebra. Further, the $R(T)$-module $K^T(G/B)$ is free, with
basis the $[S(x)]$  ($x\in W$); and the $R(T)$-bilinear map
$$
K^T(G/B)\times K^T(G/B)\to R(T),~(u,v)\mapsto \chi(G/B,u\cdot v)
$$
is a perfect pairing, where $u\cdot v$ denotes the product in
$K^T(G/B)$, and $\chi(G/B,-)$ denotes the equivariant Euler
characteristic; see \cite[3.39, 4.9]{KK}.

\begin{corollary}
The classes 
$$
[S^-(x)]^0=\sum_{y\in W,\;y\geq x}
(-1)^{\ell(y)-\ell(x)} [S^-(y)]\eqno(x\in W)
$$ 
form the dual basis of the basis of the $[S(x)]$ ($x\in W$).
\end{corollary}

\begin{proof} 
Observe that
$$\displaylines{
\chi(G/B,u\cdot v)=\chi(G/B\times G/B,[diag(G/B)]\cdot(u\times v))
\hfill\cr\hfill
=\sum_{x\in W}\chi(G/B,u\cdot [S(x)])\;
\chi(G/B,v \cdot [S^-(x)]^0)
\cr}$$
where the second equality follows from Corollary \ref{diag}. In
particular, we have for $y\in W$:
$$
\chi(G/B,u\cdot [S(y)])=\sum_{x\in W}\chi(G/B,u\cdot [S(x)])\;
\chi(G/B,[S(y)]\cdot [S^-(x)]^0). 
$$
As the $R(T)$-linear forms $u\mapsto \chi(G/B,u\cdot [S(x)])$ are
linearly independent, we obtain
$$
\chi(G/B,[S(y)]\cdot [S^-(x)]^0)=\delta_{x,y}.
$$
\end{proof}

As another consequence of the determination of the class of the
diagonal, we recover a formula of Mathieu for the character of the
$G$-module $H^0(G/B,\lambda+\mu)$ as a function of the dominant
weights $\lambda$ and $\mu$ \cite[Cor. 7.7]{M}. Further, we
determine the dimension of the modules $M(\lambda)$ introduced in
Section 4. 

\begin{corollary}\label{sep}
For any weights $\lambda$ and $\mu$, and for any $t\in T$, we have
$$
\chi(G/B,\lambda+\mu)(t)
=\sum_{x\in W}\chi(S(x),\lambda)(t)\;
\chi(S(xw_0),\cI_{\partial S(xw_0)}\otimes \cL_{G/B}(-w_0\mu))(t^{-1}).
$$
Therefore, $\chi(G/B,2\lambda)(t)=c_{\lambda}(t,t^{-1})$
for all $t\in T$. Further, for $\lambda\in\tcX^+$, one has 
$$
ch\;M(\lambda)(t,t^{-1})=ch\;H^0(G/B,2\lambda)(t) \text{ and } 
\dim M(\lambda)=\prod_{\alpha\in\Phi^+}
\frac{\langle 2\lambda+\rho,\check\alpha\rangle}
{\langle\rho,\check\alpha\rangle}.  
$$
\end{corollary}

\begin{proof}
Let $[\cL(\lambda,\mu)]$ be the class of the $T$-linearized line bundle
$\cL_{G/B}(\lambda)\boxtimes \cL_{G/B}(\mu)$ in 
$K^T(G/B\times G/B)$. As the restriction of this line bundle to the
diagonal is $\cL_{G/B}(\lambda+\mu)$, 
we have
$$
\chi(G/B,\lambda+\mu)=
\chi(G/B\times G/B,[diag(G/B)]\cdot[\cL(\lambda,\mu)]).
$$
By Corollary \ref{diag}, the latter is equal to
$$
\sum_{w\in W}\chi(S(w),\lambda)\;
w_0\chi(S(w_0w),\cI_{\partial S(w_0w)}\otimes\cL_{G/B}(\mu)).
$$
To complete the proof of the first equality, it suffices to check that
$$
\chi(S(w_0w),\cI_{\partial S(w_0w)}\otimes\cL_{G/B}(\mu))(w_0t)=
\chi(S(ww_0),\cI_{\partial S(ww_0)}\otimes\cL_{G/B}(-w_0\mu))(t^{-1}).
$$ 
For this, using Serre duality as in the proof of Theorem
\ref{recip}, we obtain
$$
\chi(S(w_0w),\cI_{\partial S(w_0w)}\otimes\cL_{G/B}(\mu))(w_0t)=
(-1)^{N-\ell(w)}\rho(t)\chi(S(w_0w),-\rho-\mu)(-w_0t).
$$
Further, the Demazure character formula implies that
$$
\chi(S(w_0w),\nu)(-w_0t)=\chi(S(ww_0),-w_0\nu)(t)
$$
for all weights $\nu$. It follows that
$$\displaylines{
\chi(S(w_0w),\cI_{\partial S(w_0w)}\otimes\cL_{G/B}(\mu))(w_0t)=
(-1)^{N-\ell(w)}\rho(t)\chi(S(ww_0),-\rho+w_0\mu)(t)
\hfill\cr\hfill
=\chi(S(ww_0),\cI_{\partial S(ww_0)}\otimes\cL_{G/B}(-w_0\mu))(-t)
\cr}$$
by Serre duality once more.

Now the second equality follows from formula $(*)$ in the proof of
Theorem \ref{recip}. For $\lambda\in\tcX^+$, the third equality
follows from the vanishing of the $H^i(Z,\lambda)$ (\cite[Th.2]{R1}), 
and the fourth one from Weyl's dimension formula. 
\end{proof}

\section{Large Schubert varieties are Cohen-Macaulay}

In this section, we prove the statement of the title and we give some
applications. We begin by constructing a partial desingularisation of
$\bX(w)$, by the total space of a fibration with fiber $\oB$ over the
usual Schubert variety $S(w)$.

For this, consider the action of $B$ on $\oB$ by left
multiplication, and the associated fiber bundle
$G\times_B\oB$ over $G/B$. The map
$G\times\oB\to\bX:~(g,x)\mapsto gx$ defines a 
birational, $G\times B$-equivariant morphism 
$$
\varphi:G\times_B\oB\to\bX
$$ 
where the action of $G\times B$ on $G\times_B\oB$ is defined by
$(g,b)(g',x)=(gg',xb^{-1})$. On the other hand, the projection 
$$
\psi:G\times_B\oB\to G/B
$$ 
is a locally trivial fibration with fiber $\oB$. Observe that 
$(\varphi,\psi)$, being the composition of 
$$
G\times_B\oB \hookrightarrow G\times_B \bX \cong G/B \times \bX, 
$$
is a closed embedding.

Let $\bX'(w)$ be the preimage of $S(w)$ under $\psi$; then
$\bX'(w)$ is stable by the subgroup $B\times B$ of 
$G\times B$. Observe that $\bX'(w)$ is the closure of 
$BwB\times_B B\simeq BwB$ in $G\times_B\oB$. As a
consequence, $\varphi$ restricts to a $B\times B$-equivariant morphism
$$
f:\bX'(w)\to\bX(w)
$$ 
which is an isomorphism over $BwB$. Denote by $\partial\bX(w)$ the
complement of $BwB$ in $\bX(w)$, and by $\partial\bX'(w)$ its
preimage under $f$. Finally, let 
$$
g:\bX'(w)\to S(w)
$$ 
be the restriction of $\psi$. Then $g$ is a locally trivial fibration
with fiber $\oB$, too.

\begin{theorem}\label{cm1} 
With notation as above, we have:

\noindent
(i) $\bX'(w)$ is Cohen-Macaulay with canonical
sheaf $\cI_{\partial\bX'(w)}\otimes f^*\cL_{\bX(w)}(-\rho)[\rho,\rho]$.

\noindent
(ii) $f_*\cO_{\bX'(w)}=\cO_{\bX(w)}$, 
$f_*\omega_{\bX'(w)}=\omega_{\bX(w)}$ and the higher direct
images $R^i f_*\cO_{\bX'(w)}$, $R^i f_*\omega_{\bX'(w)}$
vanish for $i\geq 1$.

\noindent
(iii) $\bX(w)$ is Cohen-Macaulay with canonical sheaf 
$\cI_{\partial\bX(w)}\otimes\cL_{\bX(w)}(-\rho)[\rho,\rho]$.

\noindent
(iv) The graded ring 
$R(w)=\oplus_{\lambda\in\tcX}\;H^0(\bX(w),\lambda)$ 
is Cohen-Macaulay.
\end{theorem}

\begin{proof}  
Because $S(w)$ and $\oB$ are Cohen-Macaulay, the same holds
for $\bX'(w)$. And because $f$ is birational and $\bX(w)$ is normal,
we have $f_*\cO_{\bX'(w)}=\cO_{\bX(w)}$.

We now show that $R^i f_*\cO_{\bX'(w)}=0$ for $i\geq 1$.
For this, it suffices, by a lemma of Kempf (see
e.g. \cite[II.14.13]{J}), to show that
$H^i(\bX'(w),f^*\cL_{\bX(w)}(\lambda))$ $=0$ for $i\geq 1$ and for any
regular dominant weight $\lambda$.
Consider the line bundle $\varphi^*\cL_{\bX}(\lambda)$ and its higher
direct images $R^j\psi_*(\varphi^*\cL_{\bX}(\lambda))$ for 
$j\geq 0$. Then $R^j\psi_*(\varphi^*\cL_{\bX}(\lambda))$ is the
$\tG$-linearized sheaf on $G/B=\tG/\tB$ associated with the
$\tB$-module $H^j(\oB,\lambda)$, and 
$R^j g_*(f^*\cL_{\bX(w)}(\lambda))$ is the 
restriction to $S(w)$ of this $\tG$-linearized sheaf. As
$H^j(\oB,\lambda)=0$ for all $j\geq 1$ by Corollary \ref{sur}, we have 
$R^j g_*(f^*\cL_{\bX(w)}(\lambda))=0$ for $j\geq 1$.

For a $\tB$-module $M$, denote by $\underline{M}$ the corresponding
homogeneous vector bundle on $G/B$. Then we obtain from the
Leray spectral sequence for $g$ that 
$$
H^i(\bX'(w),f^*\cL_{\bX(w)}(\lambda)) \cong 
H^i(S(w),g_*f^*\cL_{\bX(w)}(\lambda)) \cong 
H^i(S(w),\underline{H^0(\oB,\lambda)}) . 
$$
By Theorems \ref{fil} and \ref{vdk}, the left $\tB$-module
$H^0(\oB,\lambda)$ has a filtration with associated 
graded a direct sum of $P(\mu)$'s for certain dominant weights
$\mu$. Further, we have for $i\geq 1$: 
$$
H^i(S(w),\underline{P(\mu)})=0
$$ 
as follows from  
\cite[Prop. 1.4.2]{P} or \cite[Lemma 3.1.12]{vdK2}. Thus, 
$H^i(\bX'(w),f^*\cL_{\bX(w)}(\lambda))$ $=0$ 
and, therefore, $R^i f_*\cO_{\bX'(w)}=0$ 
for $i\geq 1$.

\smallskip 
We now determine the canonical sheaf $\omega_{\bX'(w)}$; we begin with
the relative canonical sheaf $\omega_g$ of $g:\bX'(w)\to S(w)$. 
Observe that the relative canonical sheaf of 
$\psi:G\times_B\oB\to G/B$ equals 
$\varphi^*\cL_{\bX}(-\beta-\rho)\otimes \psi^*\cL_{G/B}(\rho)[\rho]$
as a $\tG\times \tB$-linearized sheaf, where $[\rho]$ denotes the shift by
$\rho$ of the $\tB$-linearization. Indeed, $\omega_{\psi}$ is the
$(\tG\times \tB)$-linearized sheaf on $G\times_B\oB$ associated with the
$(\tB\times \tB)$-linearized sheaf $\omega_{\oB}$ on $\oB$. On the other
hand, the sheaf
$\varphi^*\cL_{\bX}(-\beta-\rho)\otimes\psi^*\cL_{G/B}(\rho)[\rho]$ 
is $\tG\times\tB$-linearized, and 
the associated $\tB\times\tB$-linearized sheaf on $\oB$ is
$\omega_{\oB}$ by Theorem \ref{gor}. 
As $g:\bX'(w)\to S(w)$ is the pull-back of $\psi$ under the
inclusion $S(w)\to G/B$, it follows that
$\omega_g=f^*\cL_{\bX(w)}(-\beta-\rho)\otimes g^*\cL_{S(w)}(\rho)[\rho]$. 
In particular, $\omega_g$ is invertible. Thus,
$$
\omega_{\bX'(w)}=g^*\omega_{S(w)}\otimes\omega_g=
g^*\cI_{\partial S(w)}\otimes f^*\cL_{\bX(w)}(-\beta-\rho)[\rho,\rho].
$$

We now claim that
$$
g^*\cI_{\partial S(w)}\otimes f^*\cL_{\bX(w)}(-\beta)
=\cI_{\partial\bX'(w)}. \eqno (1) 
$$ 
For this, observe that $\partial\bX'(w)$ contains the preimage under
$f$ of $\bX(w)\cap(\bX\setminus G)=\bX(w)\cap(D_1\cup\cdots\cup D_r)$,
a Cartier divisor. Further, the complement 
$\partial\bX'(w) \cap f^{-1}(G)$ of that divisor is
equal to $g^{-1}(\partial S(w))\cap f^{-1}(G)$. As the line bundle
associated with $\bX(w)\cap(D_1\cup\cdots\cup D_r)$ is
$\cL_{\bX(w)}(-\beta)$, it follows that
$$
\cI_{\partial\bX'(w)}=
g^*\cI_{\partial S(w)}\otimes
f^*\cI_{\bX(w)\cap(D_1\cup\cdots\cup D_r)}=
g^*\cI_{\partial S(w)}\otimes f^*\cL_{\bX(w)}(-\beta), 
$$
which proves the claim. 
We conclude that 
$$
\omega_{\bX'(w)} = \cI_{\partial\bX'(w)}\otimes f^*\cL_{\bX(w)}(-\rho)[\rho,\rho]. \eqno (2)
$$ 
 
{F}urther, since $f_*\cO_{\bX'(w)}=\cO_{\bX(w)}$ and 
$f(\partial\bX'(w)) = \partial\bX(w)$, then 
$$
f_*\cI_{\partial\bX'(w)}=\cI_{\partial\bX(w)}, \eqno (3)
$$
and, therefore, 
$$
f_*\omega_{\bX'(w)}=
\cI_{\partial\bX(w)}\otimes\cL_{\bX(w)}(-\rho)[\rho,\rho]. \eqno (4)
$$

We now prove that $R^if_*\omega_{\bX'(w)}=0$ for 
$i\geq 1$. Using Kempf's lemma, again, it suffices to prove that 
$$
H^i(\bX'(w),g^*\cI_{\partial S(w)}\otimes
f^*\cL_{\bX(w)}(\lambda-\beta))=0 \eqno (5)
$$
for $i\geq 1$ and for $\lambda\in\tcX^+$ big enough 
(we consider here $\omega_{\bX'(w)}\otimes f^*\cL(\lambda+\rho)$). 
We argue by
induction over $\ell(w)$, the case where $\ell(w)=0$ being obvious.

In the general case, choose a decomposition $w=sx$ where $s$ is a
simple reflection, and $\ell(x)=\ell(w)-1$. Let $P_s$ be 
the parabolic subgroup generated by $B$ and $s$. This defines the
variety
$$
\hat S(w):=P_s\times_B S(x)
$$ 
together with the map $\sigma : \hat S(w)\to S(w)$. Let  
$$
\hat\bX(w) = \hat S(w) \times_{S(w)} \bX'(w) = P_s \times_B \bX'(x) 
$$
with projections $\tau : \hat\bX(w)\to\bX'(w)$ and 
$q:\hat\bX(w)\to\hat S(w)$. Let $p:\hat\bX(w)\to\bX(w)$ be the
composition of $f$ and $\tau$, then $p$ is an isomorphism above
$BwB$. Further, $q$ is a locally trivial fibration with fiber
$\oB$. The $B\times B$-action on $\bX'(w)$ lifts to $\hat\bX(w)$,
where $1\times B$ acts trivially on $S(w)$ and $\hat S(w)$.

We claim that $\sigma_*\cO_{\hat S(w)}=\cO_{S(w)}$,
$\sigma_*\cI_{\partial\hat S(w)}=\cI_{\partial S(w)}$ and 
$R^i \sigma_*\cI_{\partial\hat S(w)}=0$ for $i\geq 1$. 
This follows from \cite{R1}. 
In more detail, consider a reduced expression for $x$ and let 
$\phi : V(w) \to S(w)$ denote the Bott-Samelson resolution associated 
to the corresponding reduced expression of $w = sx$. Observe that 
$\phi$ factors through $\sigma$, say $\phi = \sigma\theta$. 
Since $\phi$, $\theta$, $\sigma$ are proper and birational and $S(w)$,
$\hat S(w)$ are normal, then $\phi_*\cO_{V(w)} = \cO_{S(w)}$,  
$\theta_*\cO_{V(w)} = \cO_{\hat S(w)}$ and 
$\sigma_*\cO_{\hat S(w)} = \cO_{S(w)}$. 
Let $\partial V(w)$ denote the complement of the open $B$-orbit in
$V(w)$, then $\theta(\partial V(w)) = \partial\hat S(w)$ and 
$\sigma\theta(\partial V(w)) = \partial S(w)$, so that 
$\sigma_*\cI_{\partial\hat S(w)} =  \cI_{\partial S(w)}$. 
  
Further, by \cite[Prop. 2, Th. 4]{R1}, one has 
$$
\omega_{V(w)} \cong \cI_{\partial V(w)}\otimes \phi^*\cL_{S(w)}(-\rho), 
\quad 
\omega_{S(w)} \cong \phi_*\omega_{V(w)} \cong 
\cI_{\partial S(w)}\otimes \cL_{S(w)}(-\rho), 
$$
and 
$$
R^i \phi_*\cO_{V(w)} = 0 = R^i \phi_*\omega_{V(w)} 
$$
for $i\geq 1$. 

Since $\sigma$ is proper with fibres being points or projective lines, then 
$R^i \sigma_*\cI_{\partial\hat S(w)} = 0$ for $i\geq 2$ and, 
therefore, one obtains, by using the projection formula, that 
$$
R^1 \sigma_*(\cI_{\partial\hat S(w)}) \otimes \cL_{S(w)}(-\rho)) \cong  
R^1\sigma_*(\theta_* \omega_{V(w)}) \hookrightarrow  
R^1\phi_* \omega_{V(w)} = 0. 
$$ 
This proves the claim.

\smallskip 
Since $\tau$ is the pull-back of $\sigma$ under the locally trivial
fibration $g$ then, using again the projection formula, it follows that 

$$
\displaylines{  
R^i \tau_*(q^* \cI_{\partial\hat S(w)} \otimes
p^*\cL_{\bX(w)}(\lambda-\beta) ) 
\cong \hfill\cr\hfill \cong 
(R^i \tau_* q^* \cI_{\partial\hat S(w)}) \otimes 
f^*\cL_{\bX(w)}(\lambda-\beta)  
= \begin{cases} g^* \cI_{\partial S(w)} \otimes 
f^*\cL_{\bX(w)}(\lambda-\beta) & \text{ if } 
i = 0;\\ 0 & \text{ if } i \geq 1.\end{cases}
\cr}
$$ 
This yields  
$$
H^i(\bX'(w),g^*\cI_{\partial S(w)}\otimes f^*\cL_{\bX(w)}(\lambda-\beta)) 
\cong 
H^i(\hat\bX(w),q^*\cI_{\partial\hat S(w)}\otimes 
p^*\cL_{\bX(w)}(\lambda-\beta)), \eqno (6) 
$$
and it suffices to prove that the right-hand side vanishes  
for $i\geq 1$ and for $\lambda\in\tcX^+$ big enough.

Embed $S(x)=B\times_B S(x)$ into $\hat S(w)$, as a Cartier divisor;
then $\bX'(x)$ embeds into $\hat\bX(w)$. Observe that $\partial\hat S(w)=
S(x)\cup(P_s\times_B\partial S(x))$ whereas
$S(x)\cap(P_s\times_B\partial S(x))=\partial S(x)$. 
Thus, we have an exact sequence 
$$
0\to\cI_{\partial\hat S(w)}\to \cI_{P_s\times_B\partial S(x)}\to
\cI_{\partial S(x)}\otimes_{\cO_{\hat S(w)}}\cO_{S(x)}\to 0.
$$
Together with the induction hypothesis, it yields an exact sequence
$$
\hskip-20pt
\begin{array}{l}
H^0(\hat\bX(w),q^*\cI_{P_s\times_B\partial S(x)}\otimes 
p^*\cL_{\bX(w)}(\lambda-\beta)) \to
H^0(\bX'(x),g^*\cI_{\partial S(x)}\otimes f^*\cL_{\bX(w)}(\lambda-\beta)) \to
\vspace{4pt} \\ \phantom{H} 
\to 
H^1(\hat\bX(w),q^*\cI_{\partial\hat S(w)}\otimes 
p^*\cL_{\bX(w)}(\lambda-\beta))\to
H^1(\hat\bX(w),q^*\cI_{P_s\times_B\partial S(x)}\otimes 
p^*\cL_{\bX(w)}(\lambda-\beta))
\end{array} \eqno (7) 
$$
and isomorphisms for $i\geq 2$:
$$
H^i(\hat\bX(w),q^*\cI_{\partial\hat S(w)}\otimes 
p^*\cL_{\bX(w)}(\lambda-\beta)) \cong 
H^i(\hat\bX(w),q^*\cI_{P_s\times_B\partial S(x)}\otimes 
p^*\cL_{\bX(w)}(\lambda-\beta)).
$$
Consider the projection 
$$
\pi:\hat\bX(w)=P_s\times_B \bX'(x)\to P_s/B.
$$ 
Then the higher direct image sheaf
$R^j\pi_*(q^*\cI_{P_s\times_B\partial S(x)}\otimes 
p^*\cL_{\bX(w)}(\lambda-\beta))$
is the homogeneous vector bundle on $P_s/B$ associated with the
$B$-module 
$$
H^j(\bX'(x),g^*\cI_{\partial S(x)}\otimes f^*\cL_{\bX(w)}(\lambda-\beta)).
$$
The latter vanishes for $j\geq 1$ and large $\lambda$, by the
induction hypothesis. As $P_s/B$ is a projective line, it follows that
$$
H^i(\hat\bX(w),q^*\cI_{P_s\times_B\partial\hat S(w)}\otimes 
p^*\cL_{\bX(w)}(\lambda-\beta))=0
$$
for $i\geq 2$. And setting 
$$
M:=H^0(\bX'(x),g^*\cI_{\partial S(x)}\otimes 
f^*\cL_{\bX(w)}(\lambda-\beta)),
$$
then $(7)$ gives an exact sequence
$$
H^0(P_s/B,\underline{M})\to M\to
H^1(\hat\bX(w),q^*\cI_{\partial\hat S(w)}\otimes p^*\cL_{\bX(w)}(\lambda-\beta))
\to H^1(P_s/B,\underline{M}).
$$
To complete the proof, it remains to show that 
$H^1(\hat\bX(w),q^*\cI_{\partial\hat S(w)}\otimes p^*\cL_{\bX(w)}(\lambda-\beta))=0$. 
For this, it is enough to check that
$\underline{M}$ is generated by its global sections, that is, that 
$M$ is the quotient of a $P_s$-module. Now, using $(1)$ and $(3)$, observe that
$$
M \cong H^0(\bX(x),\cI_{\partial\bX(x)}\otimes \cL_{\bX(w)}(\lambda)).
$$
Further, $\partial\bX(x)=\bX(x)\cap P_s\partial\bX(x)$
(indeed, $\partial\bX(x)$ is obviously contained in
$P_s\partial\bX(x)\cap\bX(x)$; and $\bX(x)$ is not contained in 
$P_s\partial\bX(x)$, because $\bX(x)$ is not stable by $P_s$),
and this intersection is reduced as large Schubert varieties are
compatibly split in $\bX$. Therefore, 
$\cI_{\partial\bX(x)}=
\cI_{P_s\partial\bX(x)}\otimes_{\cO_{\bX(w)}}\cO_{\bX(x)}$,
and the restriction map
$$
H^0(\bX(w),\cI_{P_s\partial\bX(x)}\otimes\cL_{\bX(w)}(\lambda-\beta))\to
H^0(\bX(x),\cI_{\partial\bX(x)}\otimes\cL_{\bX(w)}(\lambda-\beta))=M
$$
is surjective for $\lambda$ big enough, by Serre's theorem. Thus, $M$
is a quotient of a $P_s$-module. This completes the proof of (ii). 

Now the previous arguments and Lemma \ref{dual} imply that $\bX(w)$ is
Cohen-Macaulay with canonical sheaf $f_*\omega_{\bX'(w)}$, which
proves (iii). Then (iv) follows
by arguing as in the proof of Theorem \ref{gor}.
\end{proof}

In particular, the closure in $\bX$ of any parabolic subgroup $P$ is
Cohen-Macaulay. As in Section 6, this leads to a degeneration of the
diagonal in $G/P$ into a union of Schubert varieties, and to formulae
for the class of the diagonal in $K^T(G/P\times G/P)$. 

\medskip
Consider now the subvariety $Z(w)=\bX(w)\cap Y$ of $\bX(w)$, and its
preimage $Z'(w)$ under $f:\bX'(w)\to\bX(w)$. We still denote by
$f:Z'(w)\to Z(w)$ and $g:Z'(w)\to S(w)$ the restrictions of $f$ and
$g$; then $g$ is a locally trivial fibration with fiber $Z$. 

As $Z'(w)=(\overline{BwB}\cap G)\times_B Z$, one has 
$g^{-1}\partial S(w)=\bigcup_{x<w}
(\overline{BxB}\cap G)\times_B Z$, and 
$$
f(g^{-1}\partial S(w)) = \bigcup_{x<w} \overline{BxZ}
=\bigcup_{x<w} Z(x)
=\bigcup_{\deuxind{x,y\in W}{x<w,\,\ell(xy)=\ell(x)+\ell(y)}}
S(xy)\times S(yw_0).
$$
We shall denote this subvariety of $Z(w)$ by $\delta Z(w)$. 

\begin{corollary}\label{cm2}
With notation as above, we have:

\noindent
(i) $Z'(w)$ is Cohen-Macaulay with canonical sheaf 
$g^*\cI_{\partial S(w)}\otimes f^*\cL_{Z(w)}(-\rho)[\rho,\rho]$.

\noindent
(ii) $f_*\cO_{Z'(w)}=\cO_{Z(w)}$, $f_*\omega_{Z'(w)}=\omega_{Z(w)}$
and the higher direct images $R^i f_*\cO_{Z'(w)}$, 
$R^i f_*\omega_{Z'(w)}$ vanish for $i\geq 1$.

\noindent
(iii) $Z(w)$ is Cohen-Macaulay with canonical sheaf 
$\cI_{\delta Z(w)}\otimes\cL_{Z(w)}(-\rho)[\rho,\rho]$.

\noindent
(iv) The graded ring $A(w)=\oplus_{\mu\in\tcX^+}\, H^0(Z(w),\mu)$ is
Cohen-Macaulay.
\end{corollary}

\begin{proof}  
Since $Z(w)$ is the complete intersection in $\bX(w)$
of the Cartier divisors $\bX(w)\cap D_1,\ldots,\bX(w)\cap D_r$, 
by Corollary \ref{reg}, 
it follows that $Z(w)$ is Cohen-Macaulay. Similarly, $Z'(w)$ is
Cohen-Macaulay and its canonical sheaf is the restriction to
$Z'(w)$ of
$$
\omega_{\bX'(w)}\otimes
f^*\cL_{\bX(w)}(\alpha_1)\otimes\cdots\otimes f^*\cL_{\bX(w)}(\alpha_r)=
\omega_{\bX(w')}\otimes f^*\cL_{\bX(w)}(\beta).
$$
The latter is equal to 
$g^*\cI_{\partial S(w)}\otimes f^*\cL_{\bX(w)}(-\rho)[\rho,\rho]$, as
we saw in the proof of Theorem \ref{cm1}. This proves (i).

The multiplication by $\sigma_1$ defines exact sequences
$$
0\to \cL_{\bX(w)}(-\alpha_1)\to \cO_{\bX(w)}\to
\cO_{\bX(w)\cap D_1}\to 0
$$
and 
$$
0\to f^*\cL_{\bX(w)}(-\alpha_1)\to \cO_{\bX'(w)}\to
\cO_{\bX'(w)\cap f^{-1}(D_1)}\to 0.
$$
By Theorem \ref{cm1}(ii), it follows that 
$f_*\cO_{\bX'(w)\cap f^{-1}(D_1)}=\cO_{\bX(w)\cap D_1}$ and 
$R^i f_*\cO_{\bX'(w)\cap D_1} =0$ for $i\geq 1$. Iterating this argument,
we obtain  $f_*\cO_{Z'(w)}=\cO_{Z(w)}$ and $R^i f_*\cO_{Z'(w)}=0$ for
$i\geq 1$. The vanishing of $R^i f_* \omega_{Z'(w)}$ and the equality
$f_*\omega_{Z'(w)}=\omega_{Z(w)}$ follow similarly from the exact 
sequences
$$
0\to \omega_{\bX(w)}\to\omega_{\bX(w)}\otimes\cL_{\bX(w)}(\alpha_1)
\to \omega_{\bX(w)\cap D_1}\to 0
$$
and
$$
0\to \omega_{\bX'(w)}\to\omega_{\bX'(w)}\otimes f^*\cL_{\bX(w)}(\alpha_1)
\to \omega_{\bX'(w)\cap f^{-1}(D_1)}\to 0
$$
together with Theorem \ref{cm1}(ii). This proves (ii).

It also follows, using Lemma \ref{dual}, that 
$$
\omega_{Z(w)}=f_*\omega_{Z'(w)}
=f_*g^*\cI_{\partial S(w)}\otimes\cL_{Z(w)}(-\rho)[\rho,\rho].
$$
But $g^*\cI_{\partial S(w)}=\cI_{g^{-1}(\partial S(w))}$ as $g$ is a
locally trivial fibration, and 
$f_*\cI_{g^{-1}(\partial S(w))}=\cI_{fg^{-1}(\partial S(w))}$ as
$f_*\cO_{Z'(w)}=\cO_{Z(w)}$. This completes the proof of (iii).

Finally, (iv) is checked as in the proof of Theorem \ref{gor}.
\end{proof} 

We now apply these geometric results to the structure of the
$\tB\times\tB$-modules $H^0(\bX(w),\lambda)$ and
$H^0(Z(w),\lambda)$. For this, we recall the definition of the Joseph
functors, see \cite[1.4]{P} and \cite[2.2]{vdK2}. Let $y,z\in W$ and
let $N$ (resp. $M$) be a $\tB$-module (resp. $\tB\times\tB$-module),
then
$$
H_y N:=H^0(S(y),\underline{N}) \text{ \ and \ } 
H_{y,z} M:=H^0(S(y)\times S(z),\underline{M}), 
$$
where $\underline{N}$ (resp. $\underline{M}$) is the corresponding $\tG$ 
(resp. $\tG\times \tG$) linearized vector
bundle on $G/B$ (resp. $G/B \times G/B$). Observe that $H_y M$, where
$M$ is regarded as a $\tB \times 1$-module, has a natural structure of
$\tB\times \tB$-module and, furthermore, there is an isomorphism of
$\tB\times \tB$-modules $H_y M \cong H_{y,1} M$. 

\begin{corollary}\label{ind}
For any weight $\lambda$, we have
$$
H^0(\bX(w),\lambda) \cong  H_{w,1} H^0(\oB,\lambda) \text{ \ and \ } 
H^0(Z(w),\lambda) \cong H_{w,1} M(\lambda).
$$ 
Further, each endomorphism of the $\tB\times\tB$-module
$H^0(Z(w),\lambda)$ is scalar. In particular, this module is
indecomposable.
\end{corollary}

\begin{proof} 
Recall that
$H^0(\bX(w),\lambda)=H^0(\bX(w),\cL_{\bX(w)}(\lambda))$. By Theorem
\ref{cm1}, the latter is isomorphic to 
$$\displaylines{
H^0(\bX'(w),f^*\cL_{\bX(w)}(\lambda)) \cong 
H^0(S(w),g_*f^*\cL_{\bX(w)}(\lambda))
\hfill\cr\hfill
\cong H^0(S(w),\underline{H^0(\oB,\lambda)}) \cong  
H_{w,1} H^0(\oB,\lambda).
\cr}
$$
Using Corollary \ref{cm2}, we obtain similarly that 
$H^0(Z(w),\lambda) \cong H_{w,1} M(\lambda)$. 

We prove that ${\rm End}_{\tB\times \tB}\,H^0(Z(w),\lambda)=k$ by descending
induction on $\ell(w)$. If $w=w_0$ then $Z(w)=Y$ and
$H^0(Z(w),\lambda)= H^0(G/B,\lambda)\boxtimes H^0(G/B,-w_0\lambda)$. In
this case, the assertion follows from \cite[II.2.8, II.4.7]{J}.

In the general case, let $s$ be a simple reflection such that
$\ell(sw)=\ell(w)+1$; let $\tP_s$ be the parabolic subgroup of $\tG$
generated by $\tB$ and $s$. Then, using \cite[2.2.5]{vdK2},  we obtain
that
$$
H^0(Z(sw),\lambda) \cong H_{sw} M(\lambda) \cong 
{\rm ind}_{\tB}^{\tP_s}H_w M(\lambda) \cong 
{\rm ind}_{\tB}^{\tP_s} H^0(Z(w),\lambda).
$$
Further, the natural map 
$$
{\rm ind}_{\tB}^{\tP_s}H^0(Z(w),\lambda)\to H^0(Z(w),\lambda)
$$ 
is surjective by Corollary \ref{sur}. Thus, 
${\rm End}_{\tB\times \tB}\,H^0(Z(w),\lambda)$ embeds into
$$
{\rm Hom}_{\tB\times \tB}\, 
({\rm ind}_{\tB}^{\tP_s} H^0(Z(w),\lambda),H^0(Z(w),\lambda)) \cong 
{\rm End}_{\tP_s\times \tB}({\rm ind}_{\tB}^{\tP_s}\; H^0(Z(w),\lambda)).
$$
The latter equals ${\rm End}_{\tB\times \tB}\,H^0(Z(sw),\lambda)$ by 
\cite[II.2.1.(7)]{J}, and we conclude by the induction hypothesis.
\end{proof} 

\smallskip 

\noindent
{\sc Remark.} By looking at right actions, one can also prove that 
$$
H^0(\bX(w),\lambda) \cong  H_{1,w^{-1}} H^0(\oB,\lambda)~{\rm and}~
H^0(Z(w),\lambda) \cong H_{1,w^{-1}} M(\lambda).
$$

\end{document}